









\magnification\magstephalf
\baselineskip14pt
\vsize24.0truecm 


\def\hatt{\widehat}
\def\dell{\partial}

\def\eps{\varepsilon}
\def\half{\hbox{$1\over2$}}

\def\quart{\hbox{$1\over4$}}
\def\sixth{\hbox{$1\over6$}}
\def\RR{\mathord{I\kern-.3em R}}
\def\Var{{\rm Var}}
\def\E{{\rm E}}
\def\d{{\rm d}}

\def\mtrix{\pmatrix} 
\def\midd{{\,|\,}}
\def\subsection{\medskip}
\def\sumin{\sum_{i=1}^n}
\def\de{\delta}
\def\ga{\gamma}
\def\arr{\rightarrow} 
\def\normal{{\cal N}}
\def\ref#1{{\noindent\hangafter=1\hangindent=20pt
  #1\smallskip}}          

\font\bigbf=cmbx12 
\font\csc=cmcsc10


\centerline{\bigbf Bayesian approaches} 
\centerline{\bigbf to non- and semiparametric density 
        estimation}

\smallskip 
\centerline{\bigbf [with a rejoinder to my discussants]} 

\medskip
\centerline{\bigbf Nils Lid Hjort, University of Oslo}

\smallskip
{{\smallskip\narrower\noindent\baselineskip12pt
{\csc Abstract.}  
This invited paper proposes and discusses several Bayesian attempts at 
nonparametric and semiparametric density estimation. 
The main categories of these ideas are as follows:
(1) Build a nonparametric prior around a given parametric model.
We look at cases where the nonparametric part of the 
construction is a Dirichlet process or relatives thereof. 
(2) Express the density as an additive expansion of 
orthogonal basis functions, and place priors on the coefficients. 
Here attention is given to a certain robust Hermite expansion
around the normal distribution. 
Multiplicative expansions are also considered.  
(3) Express the unknown density 
as locally being of a certain parametric form, then construct suitable 
local likelihood functions to express information content,
and place local priors on the local parameters. 

\smallskip\noindent
{\csc Key words:} \sl
Bayesian density estimation; 
Dirichlet process; 
Hermite expansions;
local likelihood; 
log-linear expansions; 
semiparametric estimation; 
smoothed Dirichlet priors 
\smallskip}} 

\bigskip
{\bf 1. Introduction and summary.}
Lindley (1972) noted in his review of general Bayesian methodology 
that Bayes\-ians up to then had been `embarrassingly silent' 
in the area of nonparametric statistics. 
Bayesian nonparametrics has enjoyed healthy progress since then, 
but the sub-fields of curve and surface estimation, 
nonparametric regression, 
semiparametric estimation problems, 
density and hazard rate estimation, 
and statistical pattern recognition,  
seem as yet to be not fully developed. 
This contrasts with the rapid growth and 
widespread routine use that can be witnessed 
in the frequentist corner of these areas.  

It is difficult to be a purist Bayesian in problems with 
many parameters, since setting the simultaneous prior is 
hard and parameter interactions can have unforeseen consequences.  
Such difficulties are even more prominent in the nonparametric case,
where the parameter space is infinite-dimensional, and the 
possibilities for construction of prior distributions are so unlimited. 
One must therefore expect a broader range of possible solutions, 
as opposed to the relatively clear-cut strategies for the parametric cases.
One must also expect difficulties on the technical level,
in that posterior calculations quickly become complicated.  
Furthermore, Diaconis and Freedman (1986a, 1986b) and others have 
given serious warnings about lurking dangers for 
nonparametric Bayesian constructions in the form of 
large-sample inconsistency, so performance properties of the resulting 
Bayes estimators, once derived, should also be investigated.  

\subsection
{\csc 1.1. The present article.} 
This paper is about non- and semiparametric density estimation. 
For recent accounts of many standard methods, 
see Scott (1992) and Wand and Jones (1994). 
Again, the vast majority of these are non-Bayesian,  
in the sense that they do not (explicitly) utilise 
any prior information about what general shapes 
or what degree of smoothness are more likely than others. 
We intend to propose and discuss 
several Bayesian approaches to the problem. 

One can perhaps argue that Bayesian methods are never quite as 
fullbloodedly nonparametric as some of the frequentist ones,
in that they after all require {\it some} prior knowledge 
and prior distributions as input. 
Some of the estimators we discuss are indeed semiparametric in nature, 
in that they in various ways build on `nonparametric uncertainty' 
around given parametric models. These estimators should have 
better performance properties than traditional nonparametric ones 
in a broad nonparametric neighbourhood around the parametric base model. 
Other methods are not geared towards any such parametric home grounds
and are therefore more naturally thought of as simply nonparametric. 
The methods to be discussed fall into three categories. 
Sections 2 and 3 treat estimators that build on Dirichlet processes 
in various forms, 
Sections 4 and 5 consider placing priors on coefficients in expansions,
while Sections 6 and 7 discuss non- and semiparametric methods 
that use locally parametric approximations. 
Some final remarks are offered in Section 8. 

\subsection
{\csc 1.2. Three groups of ideas.} 
The first such group of ideas has the Dirichlet process 
as basic building tool. 
The unknown distribution can be modelled as coming from
a straight Dirichlet process or from one of various related forms. 
The variants given attention to in Section 2 and 3 are 
smoothed Dirichlets, mixtures of Dirichlets, and pinned-down Dirichlets. 
Parts of the material of Section 2 are presumably known 
to workers in the field, but has been included 
since ready references do not seem to be available, 
since some of the later material in the paper builds 
on observations made here,
and since it is of interest to see that some of the 
Bayes solutions also pop up in the quite different framework 
of Sections 6 and 7. There is new material in Section 2.4 and Section 3,
on attempts at smoothing and pinning down the Dirichlet. 
 
The second general approach is to place priors on the 
coefficients of series expansions. 
In Section 4 we focus on additive orthogonal expansions 
for the densities themselves. These could for example be 
in terms of cosines or Legendre polynomials
in situations where the density is supported on a finite interval.
Particularly attractive from a semiparametric point of view 
are models of the form a normal density times an expansion  
in Hermite polynomials, since this allows modelling 
of uncertainty around the normal. 
In particular we discuss a special model of this sort 
which is more robust than the more immediate Hermite expansion. 
There are certain computational problems with additive expansions
of densities since the likelihood function quickly becomes 
a very large sum of products, leading us to outline a simplifying 
recursive computational scheme. In Section 5 we also consider 
additive expansions of the log-densities, that is, 
multiplicative expansions of the densities. This avoids 
some of the obstacles that face the otherwise attractive 
additive expansions of Section 4, and should also be easier
regarding computations. 

The third general class of methods we discuss, in Sections 6 and 7, 
is based on using locally parametric approximations to the true density,
and then placing priors on these local parameters. 
For a fixed $x$ we might for example view $f(t)=a\exp\{b(t-x)\}$ 
as a convenient approximation to the density for $t$ in the vicinity of $x$,
and one can place prior distributions on local level $a$ and 
local slope $b$, and perhaps even on the width of the local window
inside which the approximation is expected to be sufficiently adequate.  
The problem is to establish an adequate local likelihood function 
that makes it possible to compute the posterior distribution 
for the local parameters given the local data. 
For this we partly rely on methods recently 
developed in Hjort and Jones (1995). It is seen in Section 7 
that this general locally parametric approach leads to 
a long list of appealing special cases. 

\subsection
{\csc 1.3. Other work.} 
Bayesian density estimation means placing priors
on very large sets of distributions, and it is only to be
expected that this can be fruitfully done in many more ways 
than developed or mentioned in the present article. 
Here are some quick glances at other categories of such constructions. 

Building such priors via general P\'olya urn schemes 
was first treated systematically in Ferguson (1974).
The Dirichlet again occupies a special place. 
Some specialisations lead to continuous densities with probability 1, 
but often tiny details of the construction have too much 
influence on the posterior distribution. 
There has been recent renewed interest in some 
of the branches of P\'olya trees, 
see Mauldin, Sudderth and Williams (1992) and Lavine (1992). 
In particular Lavine shows how to 
construct a P\'olya tree with a given predictive density,
and how mixtures of them can model uncertainty around 
a parametric model. 

Mixtures of Dirichlet processes were first studied by Antoniak (1974).
Using such models along with hierarchical and 
otherwise generalised versions for density estimation 
is a current growth area; see references noted in Section 2.3. 
Modelling the logarithm of the density as a stochastic process
is done in Section 5; see references there, 
and further references in Lenk (1993). 

Maximum penalised likelihood and several similar methods,
such as splines smoothing, can be viewed as Bayesian.
See Good and Gaskins (1971, 1980) and the discussion in 
Silverman (1986, Section 5.4), for example,
in addition to remarks given in Section 2.4 below.  
The estimators discussed in Rissanen, Speed and Yu (1992) 
based on stochastic complexity also have Bayesian overtones. 
Finally I mention the Bayesian histograms of Hartigan (1995) 
of the present Valencia proceedings volume.
Since histograms are quite simplistic versions of general densities 
it follows that the ambitions of the present paper, 
if not necessarily the results, are higher. 

\bigskip
{\bf 2. Dirichlet process prior with smoothing.}
This section discusses various approaches 
based on the Dirichlet process or some of its smoothed and mixed relatives. 

\subsection
{\csc 2.1. Binned data and the Dirichlet smoothed histogram.}
Divide the interval where data fall into $k$ cells $C_1,\ldots,C_k$,
and let $N_j$ be the number of data points falling in $C_j$. 
These form a multinomially distributed vector with parameters
$(p_1,\ldots,p_k)$, where $p_j=F(C_j)$. 
Suppose $(p_1,\ldots,p_k)$ is given a Dirichlet prior distribution
with parameters $(ap_{0,1},\ldots,ap_{0,k})$, 
where $p_{0,1},\ldots,p_{0,k}$ are `prior guesses' for 
the $k$ probabilities and $a$ is the `strength of belief' parameter:
$p_j$ has mean $p_{0,j}$ and variance $p_{0,j}(1-p_{0,j})/(a+1)$. 
The posterior is the easily updated Dirichlet 
$(ap_{0,1}+N_1,\ldots,ap_{0,k}+N_k)$. 
If the underlying probability density is viewed as approximately 
constant over the $C_j$ interval, with length say $h_j$, 
then the Bayes estimate is 
$$\hatt f(x)=\E\Bigl\{{p_j\over h_j}\midd{\rm data}\Bigr\}
        ={1\over h_j}{ap_{0,j}+N_j\over a+n}
        =w_n{p_{0,j}\over h_j}+(1-w_n){N_j\over nh_j}, 
        \quad x\in C_j, \eqno(2.1)$$
where $w_n=a/(a+n)$.

Equation (2.1) is our first and simplest Bayesian density estimate,
and its structure is typical also for more advanced methods 
to come. It is a convex combination of the prior guess 
of the density and the histogram estimate $N_j/(nh_j)$,
with weights respectively $w_n=a/(a+n)$ and $1-w_n=n/(a+n)$. 
The estimate can perhaps be considered parametric or semiparametric 
or nonparametric, depending on the fine-ness of the binning, 
that is, the number of cells compared to the number of data points. 
A smoother estimate than the simple histogram-type version above 
emerges when a cell is placed symmetrically around the temporarily fixed $x$.
With such a moving cell $C(x)=[x-\half h,x+\half h]$ the result is 
$$\hatt f(x)=w_nh^{-1}\int_{C(x)}f_0(t)\,\d t
        +(1-w_n)f_n(x), \eqno(2.2)$$
say, in terms of a prior guess density $f_0$ and where 
$f_n(x)=n^{-1}\sumin h^{-1}I\{|x_i-x|\le \half h\}$ is the 
kernel estimator with a uniform kernel. 
And if the cell around $x$ is determined dynamically 
by the requirement that it should contain at least $r$ data points,
then $f_n$ is a $r$-nearest-neighbour estimate.

\subsection 
{\csc 2.2. The Dirichlet process prior.}
Ferguson (1973, 1974) introduced the Dirichlet process,
which in the present context allows one to carry out analysis like
in the previous subsection more easily and more generally, 
without having to discretise the sample space into cells. 
Let the distribution $F$ which governs the data points 
have such a Dirichlet process prior with parameter $aF_0$, 
where $F_0$ is a fixed distribution and $a$ is positive. 
The definition is that for each partition 
$B_1\cup\cdots\cup B_k$ of the sample space,  
$(F(B_1),\ldots,F(B_k))$ is Dirichlet with parameters 
$(aF_0(B_1),\ldots,aF_0(B_k))$. A basic result is that 
$F$ given the data is still a Dirichlet with updated parameter
$aF_0+\sumin\delta(x_i)=aF_0+nF_n$, where $\delta(x_i)$ is 
unit point mass at $x_i$, and $F_n$ is the empirical distribution function. 

This can be used to find a natural Bayesian density estimator. 
Consider $\bar f(x)=h^{-1}F[x-\half h,x+\half h]$, to be thought of,
for small $h$, as an approximation to the density at $x$. 
Its posterior distribution is given by the result quoted, and 
the Bayes estimate is found to be exactly as in (2.2),
with $f_0$ being the density of $F_0$.  
A smoother version of this argument is to 
use $\bar f(x)=\int K_h(t-x)\,\d F(t)$ instead,
where $K_h(z)=h^{-1}K(h^{-1}z)$ and $K$ is a given probability 
density symmetric around zero, 
referred to as a kernel function. Its posterior mean becomes 
$$\hatt f(x)
=\int K_h(t-x){a\,\d F_0(t)+n\,\d F_n(t)\over a+n}
        =w_n(f_0*K_h)(x)+(1-w_n)f_n(x), \eqno(2.3)$$
where $f_n(x)$ now signifies the more general $n^{-1}\sumin K_h(x_i-x)$. 
This is the classical nonparametric density estimator. 
The simpler version (2.2) corresponds to a uniform kernel on $[-\half,\half]$.   
The parameter $a$ of the prior is ideally set by the practising
statistician, in collaboration with the experts of the relevant 
field of application. The form of the posterior, and of the 
Bayes estimates derived above, suggest that $a$ has interpretation as 
`prior sample size' or strength of belief in the prior. 
The likelihood function for $a$, based on observed data,
can be derived, but it leads to a quite artificial estimator
due to special features of the unconditional distribution of data 
sampled from a Dirichlet process. The exact number of distinct 
data points, $D_n$, is a sufficient statistic for $a$,  
and one can prove that the maximum likelihood estimator 
is asymptotically equivalent to $D_n/\log n$. 
Results like this are more helpful in certain hierarchical 
constructions. Some data-based empirical Bayesian methods
for setting a value of $a$ are briefly discussed in Hjort (1991b). 

How should the smoothing parameter $h$ be chosen? 
This is the topic of hundreds of non-Bayesian papers in the literature.
It is not obvious how one should set up a Bayesian criterion 
for the selection of this parameter, which in the present context 
at least is an algorithmic parameter of an estimation method 
rather than a statistical parameter of a model. 
Note that $h$ plays a role both for both terms in (2.3).
It should not be too large since $\bar f(x)$ otherwise 
is too far away from the real parameter of interest;
in general, if $F$ has a smooth density $f$,
then $\bar f(x)$ above is equal to 
$\int K_h(t-x)f(t)\,\d t\simeq f(x)+\half\sigma_K^2h^2f''(x)$, 
where $\sigma_K^2$ is the variance of the kernel. 
In particular the prior guess used is about equal to 
$f_0+\half\sigma_K^2h^2f_0''$ rather than the preset $f_0$ itself. 
Neither should $h$ be too small. An explicit formula for 
the conditional variance of $\bar f(x)$ can be worked out, and is of the form
$$\Var\{\bar f(x)\midd{\rm data}\}
        ={1\over nh}{n^2\over (n+a)(n+a+1)}
         {1\over n}\sumin h^{-1}K(h^{-1}(x_i-x))^2 
         +{\rm smaller\ terms},$$ 
and the average of $h^{-1}K(h^{-1}(x_i-x))^2$ is of stable size
as $h\arr0$, namely about $R(K)\bar f(x)$, where $R(K)=\int K^2\,\d z$. 
Thus the posterior variance is essentially 
of order $(nh)^{-1}$ and the squared bias involved is of order $h^4$. 
Based on these facts various Bayesian criteria can be put up, 
leading to preferred size of order $n^{-1/5}$ for $h$.  
This agrees with standard results from the frequentist perspective. 

The choice of the kernel $K$ is generally less crucial than 
that of $h$. Minimising the approximative posterior variance
plus the squared bias, hinted at above, 
or for that matter the approximate risk function for the estimator,  
gives a result proportional to $\{\sigma_KR(K)\}^{4/5}$, 
which is minimal for the Yepanechnikov kernel $K(z)={3\over2}(1-4z^2)_+$ 
(scaled here to have support $[-\half,\half]$). 
West (1991) starts out from a certain marginalisation 
consistency criterion and shows that 
strict adherence to this implies that $K$ is of double exponential form. 

\subsection
{\csc 2.3. Smoothed and mixed Dirichlet priors as prior.} 
Above we gave $F$ a Dirichlet process prior and then 
smoothed the posterior $F$ around a given $x$ to produce the Bayesian 
density estimate (2.3). This approach makes perfect sense,
as does the answer. Nevertheless it is perhaps disturbing 
that the density $f$ itself has not been directly modelled,
and indeed under a Dirichlet prior it does not properly exist;
the random $F$ is with probability 1 a discrete distribution 
(with infinitely many random jumps at
an infinite collection of random locations). 

This motivates another approach, which is to `smooth the prior first', 
modelling $f$ as $f(x)=\int K_h(t-x)\,\d G(t)$ for a 
Dirichlet process $G$, say with parameter $aG_0$. 
This assures a well-defined random and continuous density,
if only $K$ is continuous. 
Such a density can also be represented as a countably infinite mixture,
as per the remark above; see Ferguson (1983). 
The posterior distribution for $f$, 
given a set of observations coming from this $f$, 
has been worked out and characterised 
via a mixture of Dirichlet processes by Lo (1984) and by Ferguson (1983). 
The exact posterior mean is however an enormous sum 
over all possible partitions of the data set, and its computation 
accordingly quite difficult for all but very small sample sizes. 
This problem can be dealt with, for example via a simulation-based method 
due to Kuo, see Ferguson (1983) and Kuo (1986), or via an iterative 
resampling scheme developed by Escobar, see Escobar and West (1994). 
There are still difficulties with the approach. 
The choices of $G_0$, $a$ and $h$ are problematic,
and the performance properties of the resulting estimators are 
less understood than those of (2.3). 
Further progress and more general hierarchical versions 
of these schemes are discussed in 
Florens, Mochart and Rolin (1992), West (1992),
and Escobar and West (1994). 

\def\bfp{{\bf p}}
\subsection
{\csc 2.4. Generalised Dirichlet priors.} 
In many situations there is some knowledge of the 
smoothness of the underlying distribution for data. 
This points to an inadequacy of the Dirichlet prior;
it is almost `too nonparametric' in its lack of contextual 
smoothness. Considering the discrete framework of Section 2.1 again,
for example, it is clear that the ordering of the $p_j$s 
is immaterial under a Dirichlet prior, whereas one's prior knowledge 
often would suggest neighbouring $p_j$s to be close with high probability 
(assuming the cell widths in that setting to be the same). 
Again this motivates trying to construct smoothed Dirichlet distributions,
to be used as more adequate priors in smoothing problems. 
The following is another route towards achieving this,
complementing the mixtures framework indicated above.  

Consider the following way of building a prior for 
a probability vector $\bfp=(p_1,\ldots,\allowbreak p_k)$: 
Let $(Y_1,\ldots,Y_k)$ be positive random variables, 
let $p_j=\exp(-Y_j)$, and condition on their sum being equal to 1. 
A calculation shows that the density of $(p_1,\ldots,p_{k-1})$ 
becomes proportional to 
$h(-\log p_1,\ldots,-\log p_k)\,p_1^{-1}\cdots p_k^{-1}$, 
where $p_k=1-\sum_{j=1}^{k-1}p_j$, 
in terms of the density $h(y_1,\ldots,y_k)$ of the $Y_i$s.  
As a special case of this construction, consider 
$$Y=(Y_1,\ldots,Y_k)\sim{\rm const.}\,
        \Bigl\{\prod_{i=1}^k \alpha_i\exp(-\alpha_i y_i)\Bigr\}\,
        g_0(y_1,\ldots,y_k), \eqno(2.4)$$
for a suitable positive function $g_0$. 
In this case the distribution for the $p_i$'s becomes proportional to 
$p_1^{\alpha_1-1}\cdots p_k^{\alpha_k-1}g(p_1,\ldots,p_k)$, 
where $g(p_1,\ldots,p_k)=g_0(-\log p_1,\ldots,\allowbreak -\log p_k)$.
In particular independent exponentials for the $Y_i$s,
corresponding to $g_0=g=1$, give the familiar Dirichlet distribution. 
Agree therefore to call this the 
{\it generalised Dirichlet distribution},
with parameters $\alpha_1,\ldots,\alpha_k$ and $g$, 
and write $\bfp=(p_1,\ldots,p_k)\sim{\cal GD}(\alpha_1,\ldots,\alpha_k;g)$
to indicate this. 

The idea is to use particular $g_0$ or $g$ functions
to push the Dirichlet in certain directions, so to speak.
A generally useful form is $g(\bfp)=\exp\{-\lambda\Delta(\bfp)\}$, that is, 
$$(p_1,\ldots,p_{k-1})\sim 
  {\rm const.}\,p_1^{\alpha_1-1}\cdots p_k^{\alpha_k-1}
        \exp\{-\lambda\Delta(p_1,\ldots,p_k)\}, \eqno(2.5)$$
where small values of $\Delta(p_1,\ldots,p_k)$ means that 
a certain characteristic of interest is present, and where
the penalty parameter $\lambda$ dictates the extent to which 
the characteristic is manifest; a large value of $\lambda$ 
forces realisations of $\bfp$ to have small values of $\Delta(\bfp)$.
Natural examples of such $\Delta(\bfp)$ functions include 
$\sum_{j=1}^{k-1}(p_{j+1}-p_j)^2$,   
$\sum_{j=2}^{k-1}(p_{j+1}-2p_j+p_{j-1})^2$, 
and $\sum_{j=1}^{k-1}(\log p_{j+1}-\log p_j)^2$.
Using $g$ functions with these $\Delta$ functions leads to 
`smoothed Dirichlet priors', forcing successive $p_j$s to be closer 
to each other with higher probability than under standard
Dirichlet conditions. Even a wish for unimodality can be built into
a suitable $\Delta(\cdot)$ function.  

It turns out that these generalised Dirichlet distributions 
are still conjugate priors for multinomial models. 
In fact, in the setting of Section 2.1, 
if $(p_1,\ldots,p_k)$ is ${\cal GD}(ap_{0,1},\ldots,ap_{0,k};g)$, 
then $\bfp$ given data is 
${\cal GD}(ap_{0,1}+N_1,\ldots,ap_{0,k}+N_k;g)$. 
For most choices of $g$ it is not possible to find explicit 
expressions for expected values, and if required these would 
have to be found by simulation or numerical integration. 
The mode is however reasonably easy to compute. 
The posterior mode $(p_1^*,\ldots,p_k^*)$ here, 
which is the Bayes solution under a sharp 0--1 loss function, 
is the maximiser of 
$$\sum_{j=1}^k(ap_{0,j}+N_j-1)\log p_j
        -\lambda\Delta(p_1,\ldots,p_k), \eqno(2.6)$$ 
under the sum to 1 constraint. This amounts to a further smoothing 
of the original Dirichlet-smoothed histogram $(ap_{0,j}+N_j-1)/(a+n-k)$,
making sure that $\Delta(p_1^*,\ldots,p_k^*)$ is not large. 
With the second $\Delta$-function mentioned above this would 
be quite similar to splines smoothing, for example. 

We have phrased the estimation problem in Bayesian terms,
starting with the (2.5) prior. If the $p_{0,j}$s are equal and
$a=k$, then the criterion to maximise is 
$\sum_{j=1}^kN_j\log p_j-\lambda\Delta(\bfp)$, 
which is the penalised log-likelihood. Other authors have 
used this as the starting point, wishing to smooth the 
simple maximum likelihood estimates across cells,
perhaps without being particularly Bayesian about it. 
The maximum penalised likelihood method can always 
be rephrased in Bayesian terms, as explained here. 
Relevant references in the present context of histogram smoothing
and density estimation include Good and Gaskins (1971, 1980),
Simonoff (1983), Silverman (1986, Section 5.4), and Hartigan (1995). 

Notice that when $n$ is large compared to $a$ and $\lambda$,
then the $\Delta$-term does not matter much, and the estimators 
become asymptotically equivalent to the ordinary histogram counts
$N_j/n$: the scaled differences $\sqrt{n}(p_j^*-N_j/n)$ tend to zero
in probability, regardless of $a$ and the $\lambda\Delta(\bfp)$ function
(unless $a$ or $\lambda$ is allowed to increase with $n$ 
at a $\sqrt{n}$ rate or faster). The situation is more 
delicate when the number of cells increases with $n$,
as should typically happen in the density estimation context. 
Simonoff (1983), who essentially worked with the third $\Delta(\bfp)$
function mentioned above, studied this sparse cells framework.
Results will not be given here, but we mention that 
methods developed in Hjort and Pollard (1995) are 
quite well suited to study properties of estimators defined by
maximisation of (2.6), also when the number of cells increases with $n$. 

The remarks and the construction given here relate to the discrete
framework with binned data. The point raised about lack of 
ordering information and lack of smoothness is also pertinent for
the Dirichlet process case; witness its representation as 
the normed version of a Gamma process with independent increments. 
It is therefore of interest to study continuous analogues 
of the distributions above, to be used, if possible, 
in the continuous data situations of Sections 2.2,
thereby by-passing the details of binning and so on. 
One possibility is to study limits in distributions of 
the discrete processes, say as the maximal binwidth tends to zero.  
This is not an easy problem, and some requirements must be placed 
on the $\Delta$ function in order to assure a well-defined limit process. 
The existence problem is made easier 
by the fact that the generalised Dirichlet distribution,
as defined above, is closed under combination of cells, 
in the following sense: Suppose $(p_1,\ldots,p_k)$ is 
${\cal GD}(\alpha_1,\ldots,\alpha_k;g)$,
and let $q_1$ be the sum of the first $j_1$ of the $p_j$s,
$q_2$ the sum of the next $j_2$ of the $p_j$s, and so on,
up to $q_m$, the sum of the last $j_m$ of the $p_j$s.
Then, working with (2.4) or the (2.5) version, 
$(q_1,\ldots,q_m)$ can be shown to be 
${\cal GD}(\sum_{r=1}^{j_1}\alpha_r,\sum_{r=j_1+1}^{j_1+j_2}\alpha_r,
\ldots;g^*)$, for a certain $g^*(q_1,\ldots,q_m)$ function. 
This generalised Dirichlet processes topic is not pursued here, however.

\subsection
{\csc 2.5. First semiparametric framework.}
Let there be a `background parameter' $\theta$ with some 
prior $\pi(\theta)$, and assume that $F$ for given $\theta$ 
is a Dirichlet $aF_0(\cdot,\theta)$. The $F_0$ could be a normal,
for example, and then this models nonparametric uncertainty 
around the normal model. 

For given $\theta$ it follows from previous comments that 
$F$ given data is a Dirichlet $aF_0(\cdot,\theta)+nF_n$,
and the arguments of Section 2.2 can be repeated to give a 
density estimator of type 
$$\hatt f(x,\theta)=\E\{\bar f(x)\midd{\rm data},\theta\}
=w_n\int K_h(t-x)f_0(t,\theta)\,\d t+(1-w_n)f_n(x), $$
where $f_0(\cdot,\theta)$ is the density of $F_0(\cdot,\theta)$. 
The final estimate is therefore of the form 
$$\hatt f(x)=\E\{\hatt f(x,\theta)\midd{\rm data}\}
=w_n\int K_h(t-x)\hatt f_0(t)\,\d t+(1-w_n)f_n(x), \eqno(2.7)$$
where $\hatt f_0(t)=\E\{\hatt f_0(t,\theta)\midd{\rm data}\}$
is the predictive parametric density. What needs to be found is 
the distribution of $\theta$ given data.  

To this end, for simplicity order the distinct data points as 
$x_1<\cdots<x_k$, and suppose the multiplicities are $j_1,\ldots,j_k$.
Let $A$ be the event that $X_i\in[x_1-\half\eps,x_1+\half\eps]$
for the first $j_1$ observations, that $X_i\in[x_2-\half\eps,x_2+\half\eps]$
for the next $j_2$ observations, and so on, where $\eps$ is small. Then 
$$\eqalign{
{\rm Pr}\{\theta\in[\theta_0&-\half\d\theta,\theta_0+\half\d\theta],A\}
\simeq\pi(\theta_0)\,\d\theta
\int{\rm Pr}\{A\midd F\}\,{\rm Dir}(aF_0(\cdot,\theta_0),\d F) \cr
&=\pi(\theta_0)\,\d\theta\,
\E\{F[x_1-\half\eps,x_1+\half\eps]^{j_1}\cdots 
        F[x_k-\half\eps,x_k+\half\eps]^{j_k}\midd\theta_0\} \cr
&\simeq\pi(\theta_0)\,\d \theta\,
{\Gamma(a)\over \Gamma(a+n)}
{\Gamma(af_0(x_1,\theta_0)\eps+j_1)\over \Gamma(af_0(x_1,\theta))}
\cdots{\Gamma(af_0(x_k,\theta_0)\eps+j_k)
        \over \Gamma(af_0(x_k,\theta))}, \cr}$$
by a formula for product moments in a Dirichlet distribution.
Since $\Gamma(b+j)/\Gamma(b)=b(b+1)\cdots(b+j-1)$ this shows that  
the posterior for $\theta$ is 
$$\pi(\theta\midd{\rm data})
={\rm const.}\,\pi(\theta)\prod_{\rm distinct}f_0(x_i,\theta), \eqno(2.8)$$
the product being over the distinct data values only. 
If in particular the data points are distinct, 
as in all proper continuous cases, 
then the sophisticated extra nonparametric randomness does not 
enter the result; the posterior is then exactly the same as 
the traditional one under the parametric model
(which also is the one corresponding to the $a$ parameter 
in the Dirichlet prior being equal to infinity).
This fills in the missing ingredient of the (2.7) estimator. 
Explicit formulae can be worked out for the case of a normal kernel,
a normal start family for $F_0$, and for traditionally used 
conjugate priors for $(\mu,\sigma)$.

\subsection
{\csc 2.6. Second semiparametric framework: Estimating the residual density.}
It is often useful to think of data as location plus noise,
and then modelling these terms separately. 
This is done in regression contexts, of course, but can also 
be done for a homogeneous sample. 
Let in general $X_i=T_\theta(\eps_i)$, 
where the $\eps_i$s are a sample from a common distribution $G$, say,
and $\theta$ is an unknown $p$-dimensional parameter 
in the transformation $T$.
This is taken to be a continuous increasing transformation 
for each given $\theta$ with inverse $\eps_i=T_\theta^{-1}(X_i)$. 
We think of the $\eps_i$s as residuals or normalised residuals. 
A simple example is $X_i=\mu+\sigma\eps_i$ 
and $\eps_i=(X_i-\mu)/\sigma$.  
The framework now will be to have some prior density for $\theta$
and in addition letting the nonparametric $G$ have a distribution
in the space of all distributions, centred at a suitable $G_0$, 
for example the standard normal.  
Note that $F(x)=G(T_\theta^{-1}(x))$, and the present Bayesian 
semiparametric setup has priors on both the $G$ part 
and the $T_\theta^{-1}(x)$ part.

Assume that $G$ is a Dirichlet with parameter $aG_0$. 
We are interested in $\bar f(x)=\int K_h(t-x)\,\d F(t)$,
now with $\d F(t)=\d G(T_\theta^{-1}(t))$, and can at least 
proceed as earlier for each given $\theta$, 
since then $G$ given data is a Dirichlet with parameter 
$aG_0+\sumin\delta(T_\theta^{-1}(x_i))$. 
One finds that $\E\{\bar f(x)\midd{\rm data},\theta\}$ can be written
$$\int K_h(t-x){a\,\d G_0(T_\theta^{-1}(t))+n\,\d F_n(t)\over a+n} 
        =w_n\int K_h(t-x)f_0(t,\theta)\,\d t+(1-w_n)f_n(x), $$
where $f_0(t,\theta)=g_0(T_\theta^{-1}(t))|\dell T_\theta^{-1}(t)/\dell t|$
is the parametric density of $X_i$ under the idealised $G=G_0$ conditions. 
This gives exactly the same density estimator as in (2.7), 
involving the predictive density $\E\{f_0(t,\theta)\midd{\rm data}\}$.
It further turns out that the posterior distribution of $\theta$
is exactly as in (2.8). This is really because the present model,
which is defined in a somewhat roundabout manner as far as 
the $X_i$s are concerned, {\it is} the same as in Section 2.5, 
with $F_0(t,\theta)=G_0(T_\theta^{-1}(t))$. 

One does, however, get interesting results of a different nature 
for the estimation of the residual density. 
The mean of $G(y)$ given both data and $\theta$ is a convex combination of 
$G_0(y)$ and $n^{-1}\sumin I\{T_\theta^{-1}(x_i)\le y\}$,
and the Bayes estimate is 
$$\hatt G(y)=w_nG_0(y)+(1-w_n)\,n^{-1}\sumin
        {\rm Pr}\{T_\theta^{-1}(x_i)\le y\midd{\rm data}\}. $$
The point is that this is a smooth estimate with a density,
in spite of the fact that $G$ under the stated prior model
does not have a continuous distribution. Typically,
$\theta$ given data is approximately a normal centred at the 
Bayes estimator (or for that matter the maximum likelihood estimator), 
say of the form $\normal_p\{\hatt\theta,\hatt V/n\}$. 
Thus $T_\theta^{-1}(x_i)$ given data is approximately a normal
with mean equal to the estimated residual, 
which we for typographical reasons write as 
$\hatt\eps_i=T^{-1}(x_i,\hatt\theta)$, 
and variance say $\hatt v_i^2/n$.
This leads to approximating 
${\rm Pr}\{T_\theta^{-1}(x_i)\le y\midd{\rm data}\}$ above with
$\Phi(\sqrt{n}(y-\hatt\eps_i)/\hatt v_i)$, 
and the density of $\hatt G$ becomes 
$$\hatt g(y)\simeq w_ng_0(y)+(1-w_n)\,
        n^{-1}\sumin h_i^{-1}\phi(h_i^{-1}(y-\hatt\eps_i))  
\quad {\rm where\ }h_i=\hatt v_i/\sqrt{n}. \eqno(2.9)$$
The second term uses a variable kernel density estimate 
for the estimated residuals with a normal kernel and 
$h_i=\hatt v_i/\sqrt{n}$ for the bandwidths. 
A result of similar nature is in Bunke (1987). 

It is remarkable that the present natural 
semiparametric Bayesian framework leads to such explicit advice 
for both the kernel form and in particular the bandwidths. 
In the case of $X_i=\mu+\sigma\eps_i$ with fixed $\sigma$  
the $\hatt v_i$s are equal and the estimator is a 
kernel estimate with smoothing parameter $\hatt v/\sqrt{n}$,
and this also happens to be the only allowed size 
in a treatment of West (1991), using a certain marginalisation 
coherence criterion. This amounts to smoothing somewhat less than 
the standard recommendation $O(n^{-1/5})$ that falls out from 
traditional mean squared error theory. 
For the case of $X_i=\mu+\sigma\eps_i$ with 
both parameters unknown, and a normal model as null model,  
the $h_i$ is approximately equal to 
$\{1+\half(x_i-\hatt\mu)/\hatt\sigma^2\}^{1/2}/\sqrt{n}$,
advicing more smoothing outside the central area than in the middle. 

While (2.9) used a large-sample approximation for the posterior of $\theta$,
exact calculations, for the derivative of the exact $\hatt G(y)$,  
are also available for many cases of interest,
for example when $X_i=\mu+\sigma\eps_i$, $G_0$ is the standard normal,
and $(\mu,\sigma)$ has a conjugate prior. 
We also point out that the apparatus here is easily extended to
regression contexts with covariate information, say with 
$X_i=T_\theta({{\bf z}}_i,\eps_i)$. The obvious special case is 
$X_i={{\bf b'z}}_i+\sigma\eps_i$, with normalised residuals 
coming from a Dirichlet process $G$ centred at the standard normal. 
There is a mildly unpleasant surprise for the situation with 
higher-dimensional data; see Section 8.7.  

\bigskip
{\bf 3. Pinned-down Dirichlet processes.}
There are other generalisations of the Dirichlet process 
worth studying in connection with density estimation. 
Such studies are potentially of interest since the final estimators
are often useful, of course, but also since models using 
these kind of infinite-dimensional priors serve as test-beds 
for general Bayesian methodology. 
Doss (1985a, 1985b),
Diaconis and Freedman (1986a, 1986b), 
Hjort (1986, 1987, 1991a) and others
have shown that some of these constructions lead to estimators 
that are inconsistent, to mention one aspect of importance;
see also Section 8.7 below. 
In the present section we look at density estimation with
certain pinned-down Dirichlet priors, for the straight 
distribution of data themselves or for the residuals. 
Still other situations worth exploring, but not pursued here, 
include the invariance-constrained Dirichlet priors of Dalal (1979). 

\subsection
{\csc 3.1. Density estimate with a pinned-down Dirichlet.}
Let $F$ be a Dirichlet process $aF_0$, conditioned on 
having fixed values $F(B_j)=z_j$ on certain 
{\it control sets} $B_1,\ldots,B_k$, where these form a partition 
of the sample space. This typically amounts to having preset values for
certain quantiles. One can prove that this pinning down of $F$ 
amounts to splitting the Dirichlet process into $k$ 
different independent Dirichlet processes, $F=z_jF_j$ on the set $B_j$,
where $F_j$ is Dirichlet with parameter $(az_j)(F_0/z_j)$; see Hjort (1986). 
Furthermore, $F_j$ given the full data set has the same distribution
as $F_j$ given only the data values that fall in $B_j$,
that is, a Dirichlet with parameter 
$aF_0+\sumin\delta(x_i)I\{x_i\in B_j\}$. 

This makes it easy to compute the mean of $F$ given data. 
If $A$ lies within $B_j$, then 
$$\E\{F(A)\midd{\rm data}\}=z_j{aF_0(A)+\sumin I\{x_i\in A\}
        \over az_j+\sumin I\{x_i\in B_j\}}
=z_j{aF_0(A)+nF_n(A)\over az_j+nF_n(B_j)}. $$ 
If $x$ is an inner point of $B_j$, therefore, the Bayes
estimate of the $\bar f(x)=\int K_h(t-x)\,\d F(t)$ parameter,
which is close to the density, is 
$$\hatt f(x)=z_j{a\int K(z)f_0(x+hz)\,\d z+nf_n(x)
        \over az_j+nF_n(B_j)}. \eqno(3.1)$$ 
If $a$ is small compared to $n$ this is close to
$\{z_j/F_n(B_j)\}f_n(x)$, which is the traditional density estimate
times a correction to account for the known location of 
a set of quantiles. 

\subsection
{\csc 3.2. Semiparametric model with pinned-down Dirichlet.}
In the framework of Section 2.6, let there be a prior $\pi(\theta)$ 
for the parameter and suppose the residual distribution $G$ for $\eps_i$ 
is given a Dirichlet process prior $aG_0$, but pinned down to have 
$G(B_j)=z_j$ on certain control sets $B_1,\ldots,B_k$, as in Section 3.1. 
Again, $f_0(x,\theta)$ is the parametric density for $X$s 
under the ideal $G=G_0$ condition. 
As a simple example of this setup, 
envisage $X_i$ as $\mu+\sigma\eps_i$ where the normalised residuals 
$\eps_i$s come from a distribution $G$ with probability mass 
$0.90$ on $[-1.645,1.645]$ and $0.05$ on 
each of $(-\infty,-1.645)$ and $(1.645,\infty)$
(1.645 being the familiar upper 5\% point of the standard normal).
Estimating parameters in this model, 
rather than in the unconstrained model that makes no such
restriction on $G$, aims at having approximately 90\% of 
future data points in 
$[\hatt\mu-1.645\,\hatt\sigma,\hatt\mu+1.645\,\hatt\sigma]$. 
Thus the general control sets apparatus is useful in connection 
with predictive analysis. 
 
Finding the posterior density for the parameters involves 
calculations that become more complicated than those 
of Sections 2.5--2.6. The result was worked out in Hjort (1986), 
and is of the form 
$$\pi(\theta\midd{\rm data})
={\rm const.}\,\pi(\theta)L_n(\theta)M_n(\theta), \eqno(3.2)$$
where 
$$L_n(\theta)=\prod_{\rm distinct}f_0(x_i,\theta)
\quad {\rm and} \quad
  M_n(\theta)=\prod_{j=1}^m{z_j^{C_j(\theta)}
        \over \Gamma(az_j+C_j(\theta))}, \eqno(3.3)$$
writing $C_j(\theta)=nF_n(T_\theta(B_j))$ for the number of 
$\eps_i=T_\theta^{-1}(x_i)$ that fall in $B_j$. 
A simple special case is the model where $G$ is taken to
have median zero, with control sets $(-\infty,0]$ and $(0,\infty)$
and $z_1=z_2=\half$. This, with known $\sigma$, 
is the situation discussed in Doss (1985a, 1985b). 
Thus (3.2)--(3.3) is a broad generalisation of the posterior
distribution found in Doss (1985a). 

The posterior density (3.2) is unusual and interesting on 
several accounts. The $L_n(\theta)$ term is largest 
around the maximum likelihood estimators, 
while $M_n(\theta)$ is largest in areas where 
$C_j(\theta)$ is close to $z_j$ for each $j$. 
These are sometimes conflicting aims, and the Bayes estimators 
in effect try to push the maximum likelihood estimates 
so as to better achieve the $z_1,\ldots,z_m$ balance. 
The unusual feature is that the data do not wash out the prior
when $n$ grows; the $M_n(\theta)$ term is about equal in strength 
to $L_n(\theta)$; see Hjort (1986, 1987) for further results and discussion.  

Let us first concentrate on estimating the density of the $X$s. 
Let $A=[x-\half\eta,x+\half\eta]$ be a short interval containing 
a given $x$. If $\theta$ is such that $T_\theta^{-1}(x)\in B_j$,
then $F(A)=G(T_\theta^{-1}(A))=z_jG_j(T_\theta^{-1}(A))$, and 
$$\E\{F(A)\midd{\rm data},\theta\}
=z_j{aG_0(T_\theta^{-1}(A))+\sumin I\{T_\theta^{-1}(x_i)\in T_\theta^{-1}(A)\}
        \over az_j+\sumin I\{T_\theta^{-1}(x_i)\in B_j\}}. $$
Hence 
$$\E\{F(A)\midd{\rm data}\}
\simeq\sum_{j=1}^m z_j\int_{D_j(x)}
{af_0(x,\theta)\eta+nF_n(A)
\over az_j+nF_n(T_\theta(B_j))}\pi(\theta\midd{\rm data})\,\d\theta, $$
where $D_j(x)$ is the set of $\theta$ for which $T_\theta^{-1}(x)\in B_j$. 
The Bayes estimate of $\bar f(x)=\int K_h(t-x)\,\d F(t)$ is therefore 
at least close to 
$$\hatt f(x)=\sum_{j=1}^m z_j\int_{D_j(x)}
{a\bar f_0(x,\theta)+nf_n(x)\over az_j+nF_n(T_\theta(B_j))}
        \pi(\theta\midd{\rm data})\,\d\theta, $$
where $\bar f_0(x,\theta)=\int K_h(t-x)f_0(t,\theta)\,\d t$. 
In cases where $T^{-1}(x,\hatt\theta)$ lies safely inside a $B_j$ set,
and the posterior distribution is not too spread out, then 
$\pi(\theta\midd{\rm data})$ gives most of its probability mass to 
a single $D_j(x)$, and a further approximation is 
$$\hatt f(x)\simeq z_j{a\E\{\bar f_0(x,\theta)\midd{\rm data}\}+nf_n(x)
        \over az_j+nG_n(B_j)}, \quad x\in T(B_j,\hatt\theta), $$
where $nG_n(B_j)$ is the number of estimated residuals in $B_j$. 

Next we pay attention to the problem of estimating the residual density,
as in Section 2.6. For subsets $A$ of $B_j$, 
$$\E\{G(A)\midd{\rm data},\theta\}
=z_j{aG_0(A)+\sumin I\{T_\theta^{-1}(x_i)\in A\}
        \over az_j+\sumin I\{T_\theta^{-1}(x_i)\in B_j\}}
=z_j{aG_0(A)+nF_n(T_\theta(A))\over az_j+nF_n(T_\theta(B_j))}. $$
This must then be averaged with respect to the (3.2) distribution. 
By results of Hjort (1987), $\theta$ is approximately a normal
centred at the Bayes estimate $\hatt\theta$ and with a 
certain variance matrix of form $\hatt V/n$. 
Hence $T_\theta^{-1}(x_i)$, given data, can be represented as 
being approximately equal to $\hatt\eps_i+\hatt v_iN/\sqrt{n}$,
where $N$ is a standard normal. 
This is similar in structure to what we saw in Section 2.6,
but the present posterior distribution is much more complicated 
than there, as is the calculation of the Bayes estimate and the 
estimated residuals. Nevertheless, for $y$ an inner point of $B_j$,    
$$\eqalign{
\E\{G[y-\half\eta,y+\half\eta]\midd{\rm data}\} 
&\simeq z_j\,\E\Bigl\{{ag_0(y)\eta
        +\sumin I\{T_\theta^{-1}(x_i)\in y\pm\half\eta\}
\over az_j+\sumin I\{T_\theta^{-1}(x_i)\in B_j\}}\,\Big|\,{\rm data}\Bigr\} \cr
&\simeq z_j\,\E\Bigl\{{ag_0(y)\eta
        +\sumin I\{\hatt\eps_i+\hatt v_iN/\sqrt{n}\in y\pm\half\eta\}
\over az_j+\sumin I\{\hatt\eps_i+\hatt v_iN/\sqrt{n}\in B_j\}}\Bigr\} \cr}$$
for small values of $\eta$. 
The random summands appearing in the denominator are mostly 
equal to 1 with high probability, if $\hatt\eps_i$ is in $B_j$, 
or equal to 0 with high probability, if $\hatt\eps_i$ is outside $B_j$. 
An approximation to the Bayes estimate of the density for 
the residuals is accordingly 
$$\hatt g(y)=z_j\,{ag_0(y)+ng_n(y)\over az_j+nG_n(B_j)}, \eqno(3.4)$$
where $G_n$ is the empirical distribution for the estimated 
residuals, so that $nG_n(B_j)$ is the number of $\hatt\eps_i\in B_j$, 
and where $g_n(y)=n^{-1}\sumin h_i^{-1}\phi(h_i^{-1}(y-\hatt\eps_i))$
is a variable bandwidth kernel estimate with explicitly 
given bandwidths, $h_i=\hatt v_i/\sqrt{n}$. 
Again, this is quite similar to (2.9), but with differently defined 
estimated residuals and bandwidths, as noted above. 
Also note that the estimator is intent on having probability mass 
close to $z_j$ on $B_j$, as the prior requests. 

\subsection
{\csc 3.3. Generalisations.}
The framework above can be generalised to regression contexts 
and to multidimensional data, essentially since the (3.2)--(3.3)
result holds in such models. It is up to the statistician
to choose control sets, for example for prediction purposes. 
We also point out that the general treatment leads to 
quantile estimates of interest, one possibility being as follows. 
Suppose the $p$th quantile $F^{-1}(p)$ is to be estimated 
for a distribution which is thought to be not very far from the normal, 
for which the exact result is $\mu+\sigma c_p$, say, where $\Phi(c_p)=p$.  
Use control sets $(-\infty,c_p]$ and $(c_p,\infty)$ with 
$z_1=p$ and $z_2=1-p$, compute the posterior density using (3.2),
and in the end use the Bayes estimate $\hatt\mu+\hatt\sigma c_p$. 
This will typically be closer to the correct $F^{-1}(p)$ than
say the maximum likelihood solution. 

\bigskip
{\bf 4. Additive Hermite expansions.}
This section discusses Bayesian density estimators 
based on certain additive expansions. 
These expansions are valid for broad classes of densities,
and can as such be viewed as nonparametric or semiparametric,
depending on whether the number of terms used is allowed 
to be infinite (or very large) or moderate. 
The methods we give are valid for all the familiar 
expansions in terms of orthogonal basis functions, 
for example cosine expansions and Legendre polynomial expansions 
for densities with support on a finite interval.  
We focus here on expansions that use Hermite polynomials to 
`correct on the normal'. These also lead to frequentist density estimates 
of interest in their own right; see Fenstad and Hjort (1995).   
The present Bayesian programme is to place priors on 
the coefficients in these expansions and work out 
posterior moments.

\subsection
{\csc 4.1. The straight Hermite expansion.}
The Hermite polynomials are defined via the derivatives
of the standard normal density, 
$\phi^{(j)}(x)=(-1)^jH_j(x)\phi(x)$. Thus 
$H_0(x)=1$,
$H_1(x)=x$,
$H_2(x)=x^2-1$,
$H_3(x)=x^3-3x$,
$H_4(x)=x^4-6x^2+3$, and so on. 
They are orthogonal with respect to the normal density,
$\int H_jH_k\phi\,\d x=j!\,I\{j=k\}$. 
For an arbitrary density $f$, consider approximations of 
the form $f_m=\phi\sum_{j=0}^m(\ga_j/j!)H_j$. 
The best approximation, in the sense of minimising 
$\int\{f/\phi-\sum_{j=0}^m(\ga_j/j!)H_j\}^2\,\d x$, 
emerges when $\ga_j=\int H_jf\,\d x$. Thus
$\ga_0=1$, $\ga_1=\E_fX$, $\ga_2=\E_f(X^2-1)$, 
$\ga_3=\E_f(X^3-3X)$, and so on. 

This approximation can be expected 
to be most effective if $f$ at the outset is not too far from 
the starting point $\phi$. It therefore helps to pre-transform
to $Y=(X-\mu)/\sigma$, writing $\mu$ and $\sigma$ for mean and 
standard deviation, develop the approximation on that scale,
and back-transform. The result is  
$$f_m(x)=\phi\Bigl({x-\mu\over \sigma}\Bigr){1\over \sigma}
\Bigl\{1+\sum_{j=3}^m{\ga_j\over j!}H_j\Bigl({x-\mu\over \sigma}\Bigr)
        \Bigr\}, \eqno(4.1)$$
where $\ga_j=\E_fH_j((X_i-\mu)/\sigma)$. Note that 
$\ga_0=1$, $\ga_1=0$, $\ga_2=0$, while 
$$\ga_3=\E_fY^3, \quad 
\ga_4=\E_f(Y^4-3), \quad
\ga_5=\E_f(Y^5-10Y^3) $$
are equal to the skewness, kurtosis, pentakosis and so on.  
The normal density is the one having each $\ga_j=0$ for $j\ge3$. 
A natural Bayesian semiparametric scheme is to put
priors on $(\mu,\sigma)$ and the $\ga_3,\ga_4,\ldots$ coefficients,
and then calculate the posterior mean of (4.1). 
This amounts to Bayesian modelling of uncertainty around the 
normal density.  

This quickly becomes quite cumbersome due to the 
prohibitively large number of terms involved,
in principle $(m+1)^n$, when the likelihood product 
$f(x_1)\cdots f(x_n)$ is expanded as a sum. 
Some progress is nevertheless possible. 
The original version (Hjort, 1994a) of the present paper 
provides details on a computational method that 
in essence manages to reduce the number of terms needed 
in the calculations from the astronomical $(m+1)^n$ 
to the hopefully manageable $(n+1)\cdots(n+m+1)/(m!)$,
by a suitably engineered recursive regrouping. 

\subsection
{\csc 4.2. The robust Hermite expansion.}
There are some difficulties with the (4.1) expansion.
The $\ga_j$ coefficients are not always finite, 
and estimates, whether frequentist or Bayesian, are 
quite variable. Another and perhaps more serious disadvantage 
surfaces when one rewrites the underlying loss function that led to the 
(4.1) approximation as $\int(f-f_m)^2\phi^{-2}\,\d x$,
that is, integrated weighted squared error with weight function
proportional to $\exp\{(x-\mu)^2/\sigma^2\}$. 
This would mean caring too much about what happens outside 
the mainstream area. 

Another Hermite expansion that in several senses gives 
more robust estimation is developed in Fenstad and Hjort (1995). 
This robust parallel to (4.1) is 
$$f_m(x)=\phi\Bigl({x-\mu\over \sigma}\Bigr){1\over \sigma}
\sum_{j=0}^m\de_j
        H_j\Bigl(\sqrt{2}{x-\mu\over \sigma}\Bigr)/\sqrt{j!}\,, \eqno(4.2)$$
where 
$\de_j=\sqrt{2}\E_fH_j(\sqrt{2}(X-\mu)/\sigma)
        \exp\{-\half(X-\mu)^2/\sigma^2\}/\sqrt{j!}\,$. 
Note that the function being averaged is bounded, so all
coefficients automatically exist, and robust estimation is unproblematic.  
The normal density has $\de_0=1$ and $\de_j=0$ for $j\ge1$.
A Bayesian model for a random density is once more to place 
a prior distribution for $(\de_0,\ldots,\de_m)$,
for example having independent $\de_j$s, and perhaps having 
distributions concentrated around zero to model 
a density in the vicinity of the normal curve. 
The techniques mentioned in the previous subsection can now be applied
with the necessary modifications to give a method for 
computing the Bayes estimate, the posterior mean of (4.2). 

Space reasons preclude giving details here; 
such are however provided in the fuller report Hjort (1994a). 
A particular construction is given there with priors formed 
by taking independent smooth priors centred around zero 
on the $\de_j$ parameters, and with prior distribution parameters
seeming reasonable in view of experience with the (4.2) expansions 
gathered in Fenstad and Hjort (1995). The result is a method 
that for such semiparametric priors gives a curve of the form 
$$\hatt f_m(x\midd\mu,\sigma)
=\phi\Bigl({x-\mu\over \sigma}\Bigr){1\over \sigma}
\sum_{j=0}^m\hatt\de_j(\mu,\sigma)
        H_j\Bigl(\sqrt{2}{x-\mu\over \sigma}\Bigr)/\sqrt{j!}$$
for each $(\mu,\sigma)$. The final estimator emerges by averaging 
this curve over say 100 realisations of $(\mu,\sigma)$ drawn from
a separately obtained posterior distribution for these. 
One possibility for this particular ingredient of the scheme is 
$$\mtrix{\mu/\sigma^* \cr \log\sigma \cr}
\approx\normal_2\{\mtrix{\mu^*/\sigma^* \cr \log\sigma^*},
{1\over n}\mtrix{1 & \half\ga_3^* \cr
        \half\ga_3^* & \half+\quart\ga_4^* \cr}\}, \eqno(4.3)$$
in which $\mu^*$, $\sigma^*$, $\ga_3^*$ and $\ga_4^*$ are 
the usual (frequentist) estimates for mean, standard deviation,
skewness and kurtosis. This is based on the approximate
binormal sampling distribution for $(\mu^*,\sigma^*)$ given 
$(\mu,\sigma)$, and (4.3) results from starting with a 
flat reference prior for $(\mu,\log\sigma)$.
It is still valid as an approximation 
with any other continuous positive prior, if $n$ is large. 
The traditionally obtained posterior in a normal model, 
based on a normal-Gamma start for $(\mu,1/\sigma^2)$, 
is like in (4.3) above, but with skewness and kurtosis parameters 
equal to zero, since then the normal-model likelihood has been used.
Thus (4.3) is the robust modification of the traditional posterior. 

\subsection
{\csc 4.3. Remarks.} 
The methods described above are quite general in nature, 
of course, and do not rely on the specifics of the Hermite expansion. 
Another attractive framework, for example, 
this time for densities on $[0,1]$, uses 
$f_m(x)=1+\sum_{j=1}^m c_j\sqrt{2}\cos(j\pi x)$  
with $c_j=\E_f\sqrt{2}\cos(\pi jX)$, 
and with individual or simultaneous priors 
placed on the collection of $c_j$s. 
(In this connection see also Brunk and Jones (1990) 
for an approach related to ours.) 
It should also be realised that as long as the order $m$ 
is fixed the underlying orthogonality structure does not matter much either, 
as far as the mathematical and computational details are concerned;
the Hermite expansion (4.1) is for example nothing but an 
$m$th order polynomial in $(x-\mu)/\sigma$ time the normal, 
and can be represented as such. 
Nevertheless the orthogonal structure is appealing, 
since we get explicit representations of the coefficients 
in terms of the underlying density; they would otherwise change
in value and interpretation when going from order $m$ to order $m+1$, say.
This representation in terms of coordinates for orthogonal 
basis functions also invites the user to think in terms of
independent prior distributions for each coefficient,
or perhaps a joint prior with a simple dependence structure.  
This translation of prior knowledge into separate priors for coordinates  
would sometimes be too forced and not relevant enough, of course,
but it greatly simplifies the task as well as the resulting 
mathematical structure. --- Further relevant comments, 
about potential difficulties with 
the approach as well as on alternative computational methods,
are given in Hjort (1994a).  

\bigskip
{\bf 5. Log-linear expansions.}
The Bayes methods of the previous section ended up being
quite cumbersome computationally, due to the large number
of possible combinations when the likelihood was expanded. 
There are also mild problems related to the fact that the 
resulting estimates occasionally may give negative values  
and may not integrate to precisely 1. 
The present section develops some theory for 
Bayesian estimation of multiplicatively expanded densities instead.   

\def\bfc{{\bf c}}
\subsection 
{\csc 5.1. Basic framework.} 
Assume observations fall in a given bounded interval.
For a fixed set of continuous basis functions 
$\psi_1,\psi_2,\psi_3,\ldots$ on this interval, consider densities of the form 
$$f(x,\bfc)=a(\bfc)^{-1}
\exp\Bigl\{\sum_{j=1}^mc_j\psi_j(x)\Bigr\}, \eqno(5.1)$$
where $a(\bfc)=\int\exp\{\sum_{j=1}^mc_j\psi_j(x)\}\,\d x$. 
If $m$ is small this is nothing but an ordinary parametric model 
for $f$, but our intention is to let $m$ be potentially large
and even possibly infinite. The model corresponds to 
an additive expansion of the log-density in terms of the basis functions. 
Our Bayesian programme is to start out with prior distributions 
for all the $c_j$ coefficients and compute the posterior distribution
of $f(x,\bfc)$, given a sample $X_1,\ldots,X_n$. 

For a simple concrete example of such a setup,
suppose the data are scaled to fall in $[0,1]$,
and let $\psi_j(x)=x^j$ for $j\ge1$. Then (5.1) amounts to
expanding the log-density as a polynomial. 
To model uncertainty around a normal density,
say scaled so as to have at least 0.999 of its probability 
mass on $[0,1]$, we would have suitable priors on $c_1$ and $c_2$ 
(and $c_2$ would have to be negative) and in addition 
have priors concentrated around zero for $c_3,c_4,\ldots\,$. 
Note that any continuous density can be approximated 
uniformly well in this way, as a consequence of 
the Stone--Weierstrass theorem. 
This happens also in many other situations,
with a suitably engineered system of basis functions. 
We may also assume without loss of generality 
that these are uniformly bounded 
(since scale factors can be moved from $\psi_j$ to $c_j$). 
Thus (5.1) defines a bona fide density also in the infinite case 
provided $\sum_{j=1}^\infty|c_j|$ is finite. 
In the fully nonparametric case we should therefore 
make sure that this series is convergent with probability 1
under prior model circumstances. 

The comments made in Section 4.3 about representation and 
interpretation of the coefficients are valid in the present 
framework too. It aids our understanding of the parameters,
and therefore of the structure of the prior distribution to be set,
if the basis functions are made orthogonal. 
If they have been chosen to satisfy $\int\psi_j\psi_k \,\d x=I\{j=k\}$,
then $c_j=\int\psi_j(x)\log f(x)\,\d x$, and any given density
can be expanded as in (5.1) with coefficients determined from this. 
This also suggests using a prior for these where the $c_j$s 
are either independent or obey some simple dependence structure. 

\subsection
{\csc 5.2. The posterior mode.} 
We now proceed to discuss general technical matters related to 
the calculation of the Bayes estimate.    
Let $\hatt\mu_j=n^{-1}\sumin\psi_j(X_i)$ be the empirical
$\psi_j$ mean and let 
$\mu_j(c_1,\ldots,c_m)=\int\psi_j(x)f(x,\bfc)\,\d x=\E_{\bfc}\psi_j(X)$
be the theoretical mean, under the assumption of the model. 
The log-likelihood of the data is 
$n\sum_{j=1}^mc_j\hatt\mu_j-n\log a(c_1,\ldots,c_m)$, 
which is easily shown to be a concave function of the parameters. 
The maximum likelihood estimators $\hatt c_1,\ldots,\hatt c_m$ 
are the unique solutions to 
$$\hatt\mu_j=\mu_j(c_1,\ldots,c_m), \quad j=1,\ldots,m. \eqno(5.2)$$   
Now consider a prior distribution $\pi(c_1,\ldots,c_m)$ 
for the $c_j$-parameters. 
The posterior distribution is then proportional to 
$$\pi(c_1,\ldots,c_m)\exp\Bigl\{
         n\sum_{j=1}^mc_j\hatt\mu_j
        -n\log a(c_1,\ldots,c_m)\Bigr\}. \eqno(5.3)$$
In several situations it would be possible to simulate from 
this posterior distribution, thus allowing us to compute 
the conditional mean of the (5.1) density for each $x$, for example. 
In general, rather than computing this squared-error loss Bayes estimate, 
it is simpler to go for the mode of (5.3), 
giving the Bayes estimate ${\bfc}^*$ of $\bfc$ 
under a sharp 0--1 loss function, and then use $f(x,{\bfc}^*)$ in the end. 

Finding the mode of (5.3) posterior density can be cumbersome 
when $m$ is large, but one is helped by the fact that its logarithm
is often exactly or approximately concave,
which means that even simple-minded numerical optimisation 
schemes, like Newton--Raphson, will work. 
The log of (5.3) is exactly concave when $\log\pi(\bfc)$ 
is concave and otherwise approximately concave if $n$ is larger than $m$. 

\subsection
{\csc 5.3. A quadratic approximation.} 
An approximation which also sheds light on the structure of 
the solution emerges by Taylor expanding the 
$\log a(\bfc)$ term of (5.3) around 
the maximum likelihood estimator $\hatt\bfc$. One finds 
$$\eqalign{\log a(\bfc)
&\simeq\log a(\hatt\bfc)
+\sum_{j=1}^m{\dell\log a(\hatt\bfc)\over \dell c_j}(c_j-\hatt c_j) 
+\half\sum_{j,k}{\dell^2\log a(\hatt\bfc)\over \dell c_j\dell c_k}
        (c_j-\hatt c_j)(c_k-\hatt c_k) \cr
&\hskip1.65cm 
+\sixth\sum_{j,k,l}{\dell^3\log a(\hatt\bfc)\over \dell c_j\dell c_k\dell c_l}
        (c_j-\hatt c_j)(c_k-\hatt c_k)(c_l-\hatt c_l). \cr}$$
Furthermore, as a consequence of the exponential form of (5.1),  
one finds that the first, second and third order derivatives 
of $\log a(\bfc)$ appearing here are equal to respectively 
$\E_{\bfc}\psi_j(X)=\mu_j(\bfc)$, 
${\rm cov}_{\bfc}\{\psi_j(X),\psi_k(X)\}=\omega_{j,k}(\bfc)$, 
and $\E_{\bfc}\{\psi_j(X)-\mu_j(\bfc)\}\{\psi_k(X)-\mu_k(\bfc)\}
        \{\psi_l(X)-\mu_l(\bfc)\}=\gamma_{j,k,l}(\bfc)$, say. 
Disregarding additive constants, the log-posterior density 
is therefore approximately equal to 
$$\log\pi(\bfc)-\half n(\bfc-\hatt\bfc)'\Omega(\hatt\bfc)(\bfc-\hatt\bfc)
 -\sixth n\sum_{j,k,l}\gamma_{j,k,l}(\hatt\bfc)
        (c_j-\hatt c_j)(c_k-\hatt c_k)(c_l-\hatt c_l), $$
where $\Omega(\bfc)$ is the covariance matrix for the $\psi_j(X)$ variables.
To illustrate further, suppose the prior distribution for $\bfc$ 
is multinormal with mean ${\bfc}_0$ and precision matrix $\Omega_0$,
i.e.~covariance matrix $\Omega_0^{-1}$.
In cases where $n$ is moderately large, 
compared perhaps to the number of $c_j$s with significantly
spread-out prior distributions, the basic quadratic approximation
will be satisfactory and we can shave off the third order terms. 
The posterior mode ${\bfc}^*$ is then close to the minimiser of 
$\half(\bfc-{\bfc}_0)'\Omega_0(\bfc-{\bfc}_0)
        +\half n(\bfc-\hatt\bfc)'\Omega(\hatt\bfc)(\bfc-\hatt\bfc)$, 
that is, 
$${\bfc}^*\simeq\{\Omega_0+n\Omega(\hatt\bfc)\}^{-1}
        \{{\bfc}_0+n\Omega(\hatt\bfc)\hatt\bfc\}. \eqno(5.4)$$

\subsection
{\csc 5.4. Special construction.} 
Suppose $f_0(x)$ is some prior guess density and that its 
log-expansion is $\sum_{j=1}^mc_{0,j}\psi_j(x)-\log a({\bfc}_0)$.
If $m$ is allowed to be infinity and the system of basis functions 
is rich enough this expansion would be the exact $\log f_0(x)$. 
Then write 
$$f(x,\bfc)={\exp\{\sum_{j=1}^mc_{0,j}\psi_j(x)\}\over a({\bfc}_0)}
{\exp\{\sum_{j=1}^m(c_j-c_{0,j})\psi_j(x)\}
        \over a(\bfc)/a({\bfc}_0)}
=f_0(x)r(x,\bfc). $$
This rewriting in terms of a prior guess times a correction function
is conceptually helpful but not of mathematical importance for 
the final estimation method. The idea is to take $c_1,c_2,\ldots$ 
independent and place a prior on each coefficient,
as per comments made at the end of Section 5.1.  
A simple strategy is to have something like 
$c_j\sim\normal\{c_{0,j},\tau^2/j^2\}$, where $\tau$ is a 
fixed parameter determining the amount of prior uncertainty.   
The Bayesian density estimate is $f(x,{\bfc}^*)$, where 
$c_1^*,c_2^*,\ldots$ maximise 
$$\sum_{j=1}^m\{-\half j^2(c_j-c_{0,j})^2/\tau^2\}
+n\sum_{j=1}^mc_j\hatt\mu_j
-n\log\int\exp\Bigl\{\sum_{j=1}^mc_j\psi_j(x)\Bigr\}\,\d x. $$
This is again a concave function with a unique maximum. 
An approximation is provided by (5.4). Note that 
$c_j^*$ becomes close to $c_{0,j}$ for all large $j$ because of 
the increased precision $j^2/\tau^2$. 

This scheme can with some efforts be generalised to 
a more semiparametrically inspired strategy, 
with a parametrical $f_0(x,\xi)$ as initial description, 
and a nonparametric correction function $r(x,\xi,\bfc)$ to estimate.
This is not pursued here.  

We have not been very specific about the number $m$ of terms 
to include. It is neat from a puristic nonparametric point of view that 
the apparatus can handle $m=\infty$, but then the priors used 
for $c_j$s for large $j$s need to have small variances.
The observation made above, about closeness of $c_j^*$ to $c_{0,j}$,
suggests that even in cases with an infinity of priors,
stopping at a moderate $m$ could constitute 
a satisfactory numerical approximation. One natural scheme 
is to compute estimates for a numbers of $m$s, perhaps 
for $m$ up to the sample size $n$, and then select one of them
according to a suitable criterion. Such a criterion is 
the Bayesian information criterion of Schwarz (1978), 
choosing the model that maximises 
$\hbox{\csc sic}(m)=n\sum_{j=1}^m\hatt c_j\hatt\mu_j
        -n\log a(\hatt c_1,\ldots,\hatt c_m)-(\half\log n)m$. 
Here $\hatt c_1,\ldots,\hatt c_m$ are the maximum likelihood 
estimates again, but the arguments used in Schwarz' proof 
of his main result also allow these to be replaced 
by the Bayesian mode estimates $c_1^*,\ldots,c_m^*$.   
The information criterion $\hbox{\csc sic}(m)$ is really based 
on large-sample approximations, and it should also be possible to 
work out useful finite-sample modifications in the present situation. 
The $\hbox{\csc sic}$ criterion is similar to but more `stingy'
than the often used Akaike information criterion; 
Schwarz's criterion is less willing to allow many parameters. 

The normal prior specification above gives an exponent 
function which is a Gau\ss ian process. 
As such the method outlined here is related to the ones 
worked with in Lenk (1991, 1993), 
in spite of having a different starting point. 
See also Leonard (1978) and Thorburn (1986) for similarly spirited approaches. 
Lenk (1991) starts with such a Gau\ss ian process 
for the exponent and uses Karhunen--Lo\'eve representations 
to translate this into a specific form of (5.1), 
with a simple form for the simultaneous prior for the coefficients.  
The present framework and methods are in several ways
more general than in Lenk (1993), since non-normal priors 
easily can be used for the coefficients. 
The idea of building uncertainty around the normal density, 
for example, in the situation where the log-density is expanded like
a polynomial, would typically require having non-normal priors 
for at least some of the coefficients. 

\bigskip
{\bf 6. Local priors on local parametric approximations.}
The basic idea of a semiparametric approach to (non-Bayesian) 
density estimation recently developed by Hjort and Jones (1995) 
is to work with parametric models, but only trusting them locally. 
That is, a parametric vehicle model 
$f(t,\theta)$ is used for $t$ in a neighbourhood of a given $x$,
and the final density estimate is of the form $f(x,\hatt\theta(x))$,
where $\hatt\theta(x)$ is based only on data local to $x$. 
This aims at the best local parametric approximation to the true density. 
In this section semiparametric Bayesian analogues are developed,
the idea being to use local priors on the local parameters. 

\subsection
{\csc 6.1. Local likelihood for densities.}
Let $f(t,\theta)$ be a suitable parametric class of densities.
The ordinary likelihood is of course $\prod_{i=1}^n f(x_i,\theta)$.
But this conveys the information content of the data only 
when the model can be trusted fully. Suppose the density is 
modelled only locally, say $f(t)=f(t,\theta)$ for 
$t$ in the window cell $C(x)=[x-\half h,x+\half h]$. 
The simplest modification of the likelihood would be 
to only include terms with $x_i\in C(x)$, 
but this is {\it not} correct; examples illustrate that 
this is an inadequate measure of information content,
and that its maximiser is an unsatisfactory estimator of the parameter. 
A more appropriate modified likelihood is the likelihood 
that only uses information about $X_i$s to the right of $x-\half h$
and what happens to these during $[x-\half h,x+\half h]$. 
In other words, the appropriate distribution is 
the conditional one given $X_i\ge x-\half h$, 
which is $f(t,\theta)/S(x-\half h,\theta)$ for $t\in C(x)$,
in terms of the survival function $S=1-F$, 
and with probability of further surviving $x+\half h$
equal to $S(x+\half h,\theta)/S(x-\half h,\theta)$. 
Calculations provided in Hjort (1994a) 
give an exact expression for this conditional information likelihood, 
say $L_{0,n}(x,\theta)$, and satisfactory Bayesian and non-Bayesian 
methods can indeed be developed with $L_{0,n}$ as basis, 
see Hjort (1995). This machinery tends to work better with 
hazard rates than for densities, however, and it turns out to
be fruitful to work instead with a convenient approximation. 
As argued in detail in Hjort (1994a), this leads to the definition of 
the {\it local likelihood\/} for the local data in $[x-\half h,x+\half h]$, 
$$\eqalign{
L_n(x,\theta)&=\Bigl\{\prod_{x_i\in C(x)}f(x_i,\theta)\Bigr\}
        \exp\Bigl\{-n\int_{C(x)} f(t,\theta)\,\d t\Bigr\} \cr 
        &=\Bigl\{\prod_{i=1}^n 
        f(x_i,\theta)^{\bar K(h^{-1}(x_i-x))}\Bigr\}
        \exp\Bigl\{-n\int \bar K(h^{-1}(t-x))
                f(t,\theta)\,\d t\Bigr\}, \cr}\eqno(6.1)$$ 
writing for the second expression 
$\bar K(z)=1$ on $[-\half,\half]$ and zero elsewhere.
Note that the scaled version $\bar K(h^{-1}(t-x))$ 
has support $[x-\half h,x+\half h]$. 
We shall also use this with more general kernel functions $\bar K(z)$, 
requiring only that they are smooth around zero
and have `correct level' $\bar K(0)=1$. 
This also translates into $\bar K(h^{-1}(t-x))$ close to 1 
for $t$ near $x$. Typically $\bar K(z)=K(z)/K(0)$ for a 
symmetric unimodal probability density kernel function $K$. 

These matters are further discussed in Hjort (1995) and 
Hjort and Jones (1995), including other reasons for viewing 
(6.1) as carrying the local information content for data in 
the $x$-window $C(x)$, for the local parameter $\theta=\theta_x$ 
of any given parametric family $f(t,\theta)$. 
These reasons are also valid much more generally than in 
the present one-dimensional framework; 
see also Hjort (1994b) and Loader (1995). 
We call the second expression in (6.1) 
the {\it kernel smoothed local likelihood} at $x$.
The argument for including non-uniform kernels 
is that they are natural smooth generalisations 
and that the local parametric model built around $x$ is sometimes
trusted a little less a little distance from $x$ than at $x$ itself. 
The requirement about `correct level' for $\bar K$ is important.
A scale factor in $\bar K$ does not matter for estimation theory
based on maximisation of such local likelihoods, 
see Hjort and Jones (1995), 
but Bayesian consequences in what follows rest on the fact that 
if $\pi(\theta)$ is a prior for $\theta$, 
then the simultaneous density for $\theta$ and local data 
is approximately proportional to $\pi(\theta)L_n(x,\theta)$. 
Note also that when $h$ is large the local kernel smoothed likelihood 
becomes $\{\prod_{i=1}^nf(x_i,\theta)\}\exp(-n)$, 
which is the usual full likelihood 
(apart from the $\exp(-n)$ factor which is independent of the parameter). 
Thus ordinary parametric likelihood methods,
whether frequentist or Bayesian, simply correspond to the large $h$ case. 
 
Before passing to the Bayesian consequences we note 
for easy reference that
$$\eqalign{
f_n(x)&=n^{-1}\sumin h^{-1}K(h^{-1}(x_i-x))
        =K(0)(nh)^{-1}\sumin\bar K(h^{-1}(x_i-x)), \cr
g_n(x)&=n^{-1}\sumin h^{-3}(x_i-x)K(h^{-1}(x_i-x)) 
       =K(0)(nh^3)^{-1}\sumin(x_i-x)\bar K(h^{-1}(x_i-x)), \cr} \eqno(6.2)$$
for the case of $\bar K(z)=K(z)/K(0)$.
Here $f_n(x)$ is the classical kernel estimator 
already encountered in equation (2.3) and later on, 
while $g_n(x)$ similarly is a $K$-based kernel type estimator 
for the density derivative $f'(x)$ times 
the constant $\sigma_K^2$ (the variance of $K$).  
When the kernel used is $K=\phi$,
the standard normal, $g_n(x)$ is simply the derivative of $f_n(x)$.

\subsection
{\csc 6.2. Local Bayes estimates.} 
Whereas Hjort and Jones (1995) maximise these local likelihoods 
and develop theory and special cases for the resulting 
$f(x,\allowbreak\hatt\theta(x))$ estimators, the present aim is 
to develop Bayesian estimators of the form 
$$\hatt f(x)=\E\{f(x,\theta)\midd {\rm local\ data}\}. \eqno(6.3)$$
The posterior distribution in question is taken to be 
$\pi(\theta)L_n(x,\theta)/\int \pi(\theta)L_n(x,\theta)\,\d \theta$,
where $\pi(\theta)$ is the prior density for $\theta$. 
A fuller verbal description of (6.3) could be 
`a locally parametric nonparametric Bayesian density estimator',
and another informative mouthful is 
`a nonparametric Bayesian estimator of the locally best
parametric approximant to the true density'.
The approximation in question here is 
in terms of a localised form of the Kullback--Leibler distance, 
see Hjort and Jones (1995).   
When $h$ is large the local likelihood becomes the ordinary one 
and (6.3) is the familiar predictive density. When $h$ is small
the estimator is essentially nonparametric in nature. 

Note that $\theta$ changes interpretation with $x$;
for each new and temporarily fixed $x$ there is a new 
parametric approximation to $f$ near $x$ and a new prior 
for the best fitting parameter $\theta=\theta_x$. 
See also Section 8.5. 
It is also fruitful to generalise this method of local parameters
to one with both global and local parameters present, 
in the spirit of two-stage priors or hierarchical Bayesian methods. 
One such framework is to model the density as $f(x,\xi,\theta)$,
where $\xi$ is a global `background' parameter with prior $\pi_0(\xi)$ 
and $\theta$ is local to $x$. 
For each $\xi$ the method above applies and gives 
$\hatt f(x,\xi)=\E\{f(x,\xi,\theta)\midd{\rm local\ data}\}$.
The final estimator is then the average of this 
over the posterior density for $\xi$, say 
$$\hatt f(x)=\E\{\hatt f(x,\xi)\midd{\rm all\ data}\} 
        =\int\Bigl\{\int f(x,\xi,\theta)
        \pi(\theta\midd{\rm local\ data})\,\d\theta\Bigr\}
        \,\pi_0(\xi\midd{\rm all\ data})\,\d\xi. \eqno(6.4)$$
Some of the examples below are of this sort. 

\bigskip
{\bf 7. Locally parametric Bayesian density estimators: Special cases.}  
The development above, ending via the local likelihood (6.1) 
in the (6.3) and (6.4) estimators, is of course quite general,
and there is a variety of different specialisations of the method. 
Below is a partial list of interesting special cases. 
We stress again that the estimators have nonparametric 
intentions, in spite of the fact that they use local parameterisations.
When the smoothing parameter is large we are back 
in the territory of ordinary fully parametric Bayesian methods. 

\subsection
{\csc 7.1. Local constant with a Gamma prior.} 
Let the local model simply be a constant,
$f(t,\theta)=\theta$ for $t$ in a neighbourhood around $x$. 
This is unrealistic as a fine description of the density, 
but makes sense locally; the main aim is to get hold of 
the `local level' for $f$. Let furthermore 
$\theta$ have a Gamma prior $\{cf_0(x),c\}$, say,
with prior mean $f_0(x)$ and prior variance $f_0(x)/c$. 
Via (6.2) the local likelihood (6.1) is seen to take the form 
$\theta^{nhf_n(x)/k_0}\exp\{-nh\theta/k_0\}$, 
writing for simplicity $k_0=K(0)$. 
It follows that $\theta$ given the local data 
is Gamma $\{cf_0(x)+nhf_n(x)/k_0,c+nh/k_0\}$, leading to 
$$\hatt f(x)={cf_0(x)+nhf_n(x)/k_0 \over c+nh/k_0}, \eqno(7.1)$$
a weighted average of prior guess and kernel estimator. 
This is exactly as in (2.3), if the value for the prior 
strength parameter $c$ used is $ah/k_0$.  

\subsection
{\csc 7.2. Local constant with a two-stage prior.} 
To generalise, let us keep the constant level model 
$f(t,\theta)=\theta$ for $t$ near $x$, but let us employ a 
two-stage prior for $\theta$: 
$\theta$ given a certain background parameter $\xi$ 
is a Gamma $\{cf_0(x,\xi),c\}$, 
and $\xi$ has a separate background prior $\pi_0(\xi)$. 
We think of $f_0(x,\xi)$ as the background model, 
which is next to be corrected on by the data. We have 
$$\E\{f(x,\theta)\midd{\rm local\ data},\xi\}
        ={c\over c+nh/k_0}f_0(x,\xi)
        +{nh/k_0\over c+nh/k_0}f_n(x), $$
and the final Bayes estimator is 
$$\hatt f(x)=
        {c\over c+nh/k_0}\int f_0(x,\xi)\pi_0(\xi\midd{\rm data})\,\d\xi 
        +{nh/k_0\over c+nh/k_0}f_n(x). \eqno(7.2)$$
This linearly combines the predictive estimator of the 
parametric prior guess density and the nonparametric kernel estimator,
and is quite similar to the (2.7) estimator 
that was derived in a quite different framework. 
It is not difficult to find the first term explicitly 
when the background model is Gau\ss ian and the $(\mu,1/\sigma^2)$ 
is given the conjugate normal-Gamma prior, for example.    

\subsection
{\csc 7.3. Local level and local slope.} 
Another generalisation is to incorporate both local level and local slope
in the vehicle model, say 
$f(x,\theta,\beta)=\theta\exp\{\beta(t-x)\}$ for $t$ close to $x$. 
Let $\beta$ have a prior $\pi(\beta)$ reflecting prior beliefs 
about the local slope at $x$, and let 
$\theta$ be a Gamma $\{cf_0(x),c\}$ 
(where the $c$ parameter could depend on $\beta$).
The local log-likelihood can be written 
$nhf_n(x)k_0^{-1}\log\theta-nh\theta\psi(\beta h)/k_0
        +\beta nh^3g_n(x)/k_0\}$, 
in view of (6.2) again, and where 
$\psi(\beta h)=\int K(z)\exp(\beta hz)\,\d z$. Thus 
$$\theta\midd{\rm local\ data},\beta
\sim{\rm Gamma}\{cf_0(x)+nhf_n(x)/k_0,c+nh\psi(\beta h)/k_0\}, $$
and 
$$\E\{f(x,\theta,\beta)\midd{\rm local\ data},\beta\}
        ={cf_0(x)+nhf_n(x)/k_0\over c+nh\psi(\beta h)/k_0}. $$
The Bayes solution $\hatt f(x)$ is to average this over 
the posterior distribution for $\beta$ given the local data. 

The local posterior for $(\beta,\theta)$ is proportional to 
$$\pi(\beta)\,\theta^{cf_0(x)+nhf_n(x)/k_0-1}
\,\exp[-\theta\{c+nh\psi(\beta h)/k_0\}]\,\exp\{nh^3\beta g_n(x)/k_0\}. $$
Integrating out $\theta$ the result is proportional to  
$$\pi(\beta){\exp\{nh^3\beta g_n(x)/k_0\}
\over \{c+nh\psi(\beta h)/k_0\}^{cf_0(x)+nhf_n(x)/k_0}}. $$
If $K=\phi$, then $\psi(\beta h)=\exp(\half\beta^2h^2)$, 
and an approximation gives
$${\exp\{cf_n(x)\exp(-\half\beta^2h^2)\}
\over \{\exp(\half\beta^2h^2)+ck_0/(nh)\}^{cf_0(x)}}
\,\pi(\beta)\,\exp\{-\half nh^3f_n(x)\beta^2/k_0+nh^3\beta g_n(x)/k_0\}. $$
Supposing $c$ to be small compared to $nh$, and letting the prior for 
$\beta$ be a normal $(\beta_0,1/w_0^2)$, the resulting 
local posterior for $\beta$ is approximately proportional to 
$$\exp\bigl[-\half w_0^2(\beta-\beta_0)^2
-\half nh^3f_n(x)\{\beta-g_n(x)/f_n(x)\}^2/k_0\bigr], $$
that is, 
$$\beta\midd{\rm local\ data}
\approx\normal\Bigl\{{w_0^2\beta_0+nh^3f_n(x)\{g_n(x)/f_n(x)\}/k_0
        \over w_0^2+nh^3f_n(x)/k_0},
        {1\over w_0^2+nh^3f_n(x)/k_0}\Bigr\}. \eqno(7.3)$$
Note that the mean is a convex combination of $\beta_0$ and 
the natural nonparametric estimate of the log-derivative of the density.
To a first order approximation the $\beta$ given local data
is a normal with mean $g_n(x)/f_n(x)$ and variance 
$\{nh^3f_n(x)/k_0\}^{-1}$. 
The final approximation for the density estimate itself becomes 
$$\hatt f(x)\simeq f_n(x)\,
\E\exp\{-\half^2h^2\beta^2\midd{\rm local\ data}\}
\simeq f_n(x){\exp\{-\half h^2\bar\mu^2/(1+h^2\bar\sigma^2)\}
        \over (1+h^2\bar\sigma^2)^{1/2}}, $$
where $\bar\mu$ and $\bar\sigma^2$ are the newly found posterior
parameters for $\beta$. This is similar to what 
Hjort and Jones (1995, Section 5.2) found for the 
maximum local likelihood density estimate for this local log-linear model. 

\subsection
{\csc 7.4. Local level, slope and curvature.}
With efforts the previous example can be generalised 
to a model for local level, local slope and local curvature.
The calculations are as above but become more complicated,
and in addition to $f_n(x)$ and $g_n(x)$ 
they will involve an estimate of the second derivative of $f$. 
The end result is a Bayesian parallel to the maximum local likelihood 
density estimator found in Hjort and Jones (1995, Section 5.3) for
the local log-quadratic model. 
These constructions can also be generalised to 
the situation where the prior guess density involves a
background parameter, just as estimator (7.2) generalised estimator (7.1). 

\subsection
{\csc 7.5. Prior guess times a local constant.} 
This time write the density as a prior guess times 
a correction function which must then be locally estimated.
With a local constant for this purpose this means using 
$f(t,\theta)=f_0(x)\theta$ for $t$ near $x$.
Let $\theta$ have a Gamma distribution around 1,
say with parameters $(c,c)$ 
(in particular this means that $f(x,\theta)$ is seen 
as a Gamma with $\{c,c/f_0(x)\}$, and a differently 
structured variance than in 7.1 above). 
The local likelihood is 
$$\eqalign{L_n(x,\theta)
&=\prod_{i=1}^n\{f_0(x_i)\theta\}^{\bar K(h^{-1}(x_i-x))}
        \exp\Bigl\{-n\theta\int\bar K(h^{-1}(t-x))f_0(t)\,\d t\Bigr\} \cr
&=\prod_{i=1}^n\{f_0(x_i)\}^{\bar K(h^{-1}(x_i-x))}
        \,\theta^{nhf_n(x)/k_0}
        \exp\Bigl\{-nh\theta\int K(z)f_0(x+hz)\,\d z\Bigr\}, \cr}$$
which leads to a different type of posterior for $\theta$ 
than in special case 7.1, namely 
$$\theta\midd{\rm local\ data}
        \sim{\rm Gamma}\Bigl\{c+nhf_n(x)/k_0,
        c+nh\int K(z)f_0(x+hz)\,\d z/k_0\Bigr\}.$$
The density estimator becomes 
$$\hatt f(x)=f_0(x){c+nhf_n(x)/k_0
        \over c+nh\int K(z)f_0(x+hz)\,\d z/k_0}. \eqno(7.4)$$
This pushes the kernel estimator downwards in regions where $f_0$ 
is convex and upwards in regions where $f_0$ is concave. 
Again extensions are possible, to two-stage priors with 
a global parameter in $f_0(x,\xi)$, and to the local log-linear
model for $f/f_0$ with a local slope parameter in addition to $\theta$. 

\subsection
{\csc 7.6. An alternative correction factor function.} 
Again write $f(t,\theta)=f_0(t)\theta$ for a prior guess density 
$f_0$ and a local correcting constant, 
but this time consider using the alternative kernel function 
$\bar K(z)=k_0^{-1}K(z)f_0(x)/f_0(x+hz)$, that is, 
$\bar K(h^{-1}(t-x))=k_0^{-1}hK_h(t-x)f_0(x)/f_0(t)$ for $t$ near $x$. 
A calculation shows that the local kernel smoothed likelihood 
becomes proportional to 
$$\theta^{nhf_0(x)r_n(x)/k_0}
        \exp\{-nh\theta f_0(x)/k_0\}, $$
where $r_n(x)=n^{-1}\sumin K_h(x_i-x)/f_0(x_i)$ is the natural
nonparametric kernel estimator of the correction factor function 
$r(x)=f(x)/f_0(x)$. The Bayesian density estimator becomes 
$$\hatt f(x)=f_0(x){c+nhf_0(x)r_n(x)/k_0
        \over c+nhf_0(x)/k_0}. $$
This is close to $f_0(x)r_n(x)$,
which is the simplest form of a class of density estimators 
recently developed by Hjort and Glad (1995), all of the 
form initial parametric start estimate times nonparametric 
correction factor. 

An extension of the above is to write $f(t)=f(t,\xi)\theta$ 
for $t$ near $x$, where $\theta$ is local to $x$ and 
$\xi$ is a global parameter, as with (6.4). 
The $f(t,\xi)$ could for example be a normal with a prior on $(\mu,\sigma)$. 
This leads to an estimator of the form
$$\hatt f(x)=\int f(x,\xi){c+nhf(x,\xi)r_n(x,\xi)/k_0
        \over c+nhf(x,\xi)/k_0}
        \pi_0(\xi\midd{\rm data})\,\d\xi, \eqno(7.5)$$
where $r_n(x,\xi)=n^{-1}\sumin K_h(x_i-x)/f(x_i,\xi)$. 
If in particular $c$ goes to zero, 
arguably corresponding to a noninformative prior for the local constant, 
then $\hatt f(x)$ takes the form of a predictive version 
of the Hjort and Glad type estimator,
$$\int f(x,\xi)r_n(x,\xi)\pi_0(\xi\midd{\rm data})\,\d\xi
        =n^{-1}\sumin K_h(x_i-x)
        \,\E\Bigl\{{f(x,\xi)\over f(x_i,\xi)}
        \,\Big|\,{\rm data}\Bigr\}. \eqno(7.6)$$
This is simply the classical kernel estimator when the prior
model $f(x,\xi)$ is flat and noninformative, and otherwise
aims to correct the kernel estimator so as to have smaller bias
in a broad neighbourhood of the parametric model.   
The (7.5) and (7.6) estimators are Bayesian predictive versions 
of the $f(x,\hatt\xi)r_n(x,\hatt\xi)$ type estimator.
In Hjort and Glad (1995) these semiparametric estimators 
have been shown to have frequentist properties
generally comparable to and often better than 
those of the kernel estimator. 

\subsection
{\csc 7.7. A running normal density.} 
This time let us try to estimate a `running normal density'.
The local model is now a normal $(\mu,\sigma^2)$, 
and the density estimate is 
$$\hatt f(x)={1\over \sqrt{2\pi}}\int\int{1\over \sigma}
        \exp\{-\half(x-\mu)^2/\sigma^2\}\,
        \pi(\mu,\sigma\midd{\rm local\ data})\,\d\mu\,\d\sigma. \eqno(7.7)$$ 
To illustrate without too many technicalities we take $\sigma$ known.
Using the standard normal kernel a calculation shows that 
the local likelihood is proportional to 
$$\exp\Bigl[-\half{1\over \sigma^2}{nhf_n(x)\over \phi(0)}
\Bigl\{\delta-h^2{g_n(x)\over f_n(x)}\Bigr\}^2\Bigr]
\exp\Bigl\{-{nh\over \sqrt{\sigma^2+h^2}}
\exp\Bigl(-\half{\delta^2\over \sigma^2+h^2}\Bigr)\Bigr\}, $$
in terms of $\delta=\mu-x$, and using (6.2) again.  
If the local $\mu$ is given a normal $(\mu_0,\tau^2)$ prior, 
for example, then this gives the posterior density.
The (7.7) estimator is in the end found from numerical integration.

\bigskip
{\bf 8. Supplementing remarks.}
This final section gives various additional results and remarks,
several of which point to questions for further research.

\subsection
{\csc 8.1. Fine-tuning of parameters.} 
Several of the estimators arrived at in the paper depend 
on one or more parameters that need to be decided on for each application. 
Some of these parameters relate to specification of the prior 
distribution and others are `smoothing parameters'. 
Parameters in the prior should ideally be set via the 
infamous `prior considerations', 
perhaps explicitly involving previous data of similar nature, 
and sometimes they may be estimated in an empirical Bayes manner. 
Using a hyper-prior and averaging over its posterior is another 
and usually quite robust method.
 
It is also of separate interest to consider 
`noninformative' reference type priors. 
For the estimators of Sections 2 and 3 the natural version of this 
is to let the prior strength parameter $a$ tend to zero.
This leads to the traditional kernel estimator itself 
for cases (2.3) and (2.7), 
to a variable bandwidth kernel estimator using residuals in (2.9), 
to kernel estimators corrected for level in Section 3, 
and to a version of maximum penalised likelihood in (2.6). 
In particular this lends some Bayesian support to the kernel estimator
and to the other limiting versions mentioned. 
Similarly, for the schemes in Section 7 that use Gamma priors
with strength parameter $c$, the natural reference prior 
corresponds to letting $c$ tend to zero. Again this is seen
to lead to the kernel estimator for cases (7.1), (7.2), 
and to interesting competing versions for cases (7.4), (7.5), (7.6). 

The $h$ of Sections 2 and 3 is an example of an external smoothing parameter,
perhaps more appropriately seen as an algorithmic parameter
governing the amount of smoothing than a parameter of a statistical model. 
Bayesian versions of cross validation criteria can be invented 
for the purpose. 
The $h$ of Sections 6 and 7 is also such a smoothing parameter,
but it can also be interpreted in the context of a Bayesian model;  
if the kernel $K$ in question is scaled so as to have 
support $[-\half,\half]$, for example, then $[x-\half h,x+\half h]$
is the interval around $x$ in which the parametric approximation
is believed to be adequate. As such the length of this interval 
of adequacy can be given a prior, and so on. 
If a data-driven method is wished for it seems more natural, however, 
to use a suitable local goodness of fit criterion; 
the interval is stretched until the parametric model bursts,
and this defines the $h$ to be used for the given $x$. 
A specific version of this idea is described in Hjort (1994a),
using weak convergence theory to find such a bursting point. 

\subsection
{\csc 8.2. Performance.} 
With various starting constructions we have been able to 
give recipes for the computation of Bayes estimates. 
Non- and semiparametric Bayesian constructions are often so technically 
complicated that even this is sometimes an achievement. 
The question of performance analysis is typically even harder,
and is presumably the reason why this aspect of the method 
is too often neglected. 
Different Bayesians might perhaps wish to stress different aspects
of performance, but `analysing performance' in the present 
context could for example mean studying exact or approximate 
risk functions under pointwise or integrated squared error loss, 
that is, frequentist behaviour for given candidate densities. 
This might involve    
(a) assessing approximate biases and variances; 
(b) comparing these with those of standard estimators  
like the kernel method, for densities that are likely and not so
likely under the prior; 
(c) proving or disproving large-sample consistency and normality; 
(d) assessing the reduction, if any, 
in terms of Bayes risk, compared to that of standard methods, 
both under the ideal prior that led to the estimator in question
but also under other priors. 
These comments also point to simulation as a natural tool.  
   
Most estimators derived in Section 2 can be analysed like this,
at least for large $n$ and small $h$. 
These estimators are typically asymptotically equivalent
to appropriate kernel estimators,
if the prior is held fixed, and the large sample theory for
kernel estimators is very well developed. 
Estimators from Section 3 can also be analysed 
with suitable extensions of standard tools. 
Some but not all of the estimators from Section 6 and 7 
can be analysed, for moderate to large $n$ and small $h$, 
using methods of Hjort and Jones (1995). The estimators 
that rely on both global and local parameters would require 
more delicate tools for their analysis. 
%
Estimators from Sections 4 and 5 are not so easily analysed. 
The task is not terribly hard if the expansion order $m$ is fixed 
and small compared to $n$, but otherwise we find ourselves 
in need of more theory, particularly when $m$ is allowed to be infinite. 
See comments in Hjort (1994a). 
 

\subsection
{\csc 8.3. Bayesian hazard rate and regression curve estimation.} 
Many of the methods and results presented here 
have parallels in the problem of estimating hazard rates 
in survival analysis and in more general counting process
models for life history data.  
See Hjort (1991a, `third general method') and Hjort (1991b, Section 8)
for two frameworks involving nonparametric randomness around 
parametric models, the latter also extending to such Bayesian 
uncertainty around the semiparametric Cox regression model. 
In these settings, which are analogous to the present paper's
Sections 2 and 3, Beta processes, a generalised class of 
hazard function relatives of the Dirichlet (Hjort, 1990), 
play the natural role. 
A Bayesian locally parametric approach for hazard rates estimation,
analogous to Sections 6 and 7, should not be difficult to develop,
with the local likelihood machinery already present in Hjort (1995).   

The local likelihood machinery of Sections 6 and 7 
also has natural parallels in Bayesian nonparametric regression;
model the unknown regression curve as being locally linear,
place normal priors on the local line parameters, 
and compute the posterior given local data using local likelihood.
The noninformative prior versions of this scheme correspond to 
recently developed local polynomial regression methods 
that have been shown to have very good performance properties,
see for example Fan and Gijbels (1992).  
Empirical Bayesian and hierarchical Bayesian schemes 
can be developed as well, and these could easily perform better 
than standard methods in situations with several covariates. 
See Hjort (1994b).

\subsection
{\csc 8.4. Further problems.}
Standard density estimation methods work reasonably well 
in all reasonably populated areas. That is, they perform
well in dimensions 1 and possibly 2, 
but often not at all in higher dimensions 
(the curse of dimensionality is precisely that no areas
are well populated in higher dimensions), 
and often not well in the tails. 
These problems have not been focussed on in
the present paper, but are presumably areas where the Bayesian
method might improve significantly on standard methods. 
One needs a broader range of prior distributions that reflect 
various notions of what are likely and not so likely densities.
Challenges include building Bayesian methods that are geared towards 
for example approximate unimodality or approximate bimodality 
(and for this mixtures approaches would be appropriate). 
Grander problems, where the Bayesian viewpoint is well worth
exploring to a fuller extent than hitherto, 
include statistical pattern recognition,
non- and semiparametric regression (particularly with 
many covariates), and neural networks. 

\subsection
{\csc 8.5. Full Bayesian models for local parametric estimation.}
Our locally parametric method requires a local model and 
a local prior around each $x$. A full global model for the complete
density curve requires in addition a description of how these
underlying local parameters change with $x$, 
but this was not necessary as far as the computation of the
estimate is concerned. 
It is nevertheless of interest to build such a fuller 
stochastic process framework, which also would be valuable 
when it comes to estimation of prior parameters,  
and for performance evaluation, cf.~comments made in 8.1 and 8.2 above. 
This is not an easy task. To illustrate, 
take the simplest of the special cases considered, that of Section 7.1. 
Take the local constant $\theta_x$ to be the result of 
a smoothed Gamma process with independent increments,
say $\theta_x=T(x-\half\,\d x,x+\half\,\d x)/\d x$ with a certain
small smoothing resolution $\d x$, where $T(a,b)$ is Gamma distributed 
with parameters $\{c_0F_0(a,b),c_0\}$ for each given $(a,b)$ interval.
There is such a process by Kolmogorov's consistency theorem 
and the additive property of Gamma variables with equal shape parameter.    
Hence $\theta_x$ is approximately a Gamma with parameters 
$\{c_0f_0(x)\,\d x,c_0\,\d x\}$,
with mean $f_0(x)$ and variance $f_0(x)/(c_0\,\d x)$. 
This fits in with the treatment in Section 7.1. 
There are alternative constructions of interest and relevance, 
and the other special cases in Section 7 require various extensions.

\subsection
{\csc 8.6. The linear semiparametric estimator.}  
Results (2.3), (2.7) and (7.1)--(7.2) give Bayesian support to 
estimators of similar form that were considered in non-Bayesian 
frameworks by both Schuster and Yakowitz (1985) 
and Olkin and Spiegelman (1987). Jones (1994) and others have 
pointed to certain difficulties with such estimators, 
primarily because the mixing parameter is somewhat imprecisely defined 
and difficult to estimate. In the present Bayesian framework the
mixing parameter has a clear interpretation, however. 
For the case of (7.2), for example, 
a natural suggestion is to estimate $c$ by empirical Bayesian techniques
and choose $h$ by a suitable local goodness of fit criterion. 
Again methods worked out in the regression case in Hjort (1994b) 
are relevant. We leave the finer details for future work. 

\subsection
{\csc 8.7. Inconsistent estimates of residual density in higher dimensions.}
The semiparametric framework in Section 2.6 for data 
that were parametric transformations of residuals can be generalised 
to situations with multidimensional data. 
One such situation of interest is $X_i=\mu+\Sigma^{1/2}\eps_i$,
say in dimension $d$, where $\mu$ is location vector and 
$\Sigma$ is a symmetric and positive definite matrix, and 
where the $\eps_i$s come from a residual distribution $G$ in ${\RR}^d$. 
Now give $(\mu,\Sigma)$ a prior and let $G$ be a Dirichlet with
parameter $aG_0$. --- This natural setup unexpectedly 
leads to a non-consistent Bayes estimator of the residual density,
however, if $d\ge 3$. See again Hjort (1994a) for details and comments.

\bigskip
\centerline{\bf References} 

\parindent0pt
\baselineskip11pt
\parskip3pt
\medskip 

\ref{%
Antoniak, C.E. (1974).
Mixtures of Dirichlet processes with application to 
Bayesian nonparametric problems.
{\sl Annals of Statistics} {\bf 2}, 1152--1174.}

\ref{%
Brunk, H.D.~and Jones, P.W. (1990).
Fitting conditional distributions using orthogonal expansions.
{\sl Journal of Statistical Planning and Inference} 
{\bf 26}, 325--337.}

\ref{%
Bunke, O. (1987). 
Bayesian inference in semiparametric models.
Proceedings of the First World Congress of the Bernoulli Society,
Tashkent, USSR, 27--30. VNU Science Press.}  

\ref{%
Dalal, S.R. (1979).
Dirichlet invariant processes and applications to 
nonparametric estimation of symmetric distribution functions.
{\sl Stochastic Processes and Their Applications} {\bf 9}, 99--107.}


\ref{%
Diaconis, P.~and Freedman, D.A. (1986a).
On the consistency of Bayes estimates (with discussion).
{\sl Annals of Statistics} {\bf 14}, 1--67.}  

\ref{%
Diaconis, P.~and Freedman, D.A. (1986b).
On inconsistent Bayes estimates of location. 
{\sl Annals of Statistics} {\bf 14}, 68--87.}

\ref{%
Doss, H. (1985a).
Bayesian nonparametric estimation of the median.
I: Computation of the estimates.
{\sl Annals of Statistics} {\bf 13}, 1432--1444.}  

\ref{%
Doss, H. (1985b).
Bayesian nonparametric estimation of the median.
II: Asymptotic properties of the estimates. 
{\sl Annals of Statistics} {\bf 13}, 1445--1464.}

\ref{%
Escobar, M.D.~and West, M. (1994).
Bayesian density estimation and inference using mixtures.
{\sl Journal of the American Statistical Association}, to appear.} 

\ref{%
Fan, J.~and Gijbels, I. (1992).
Variable bandwidth and local linear regression smoothers.
{\sl Annals of Statistics} {\bf 20}, 2008--2036.}

\ref{%
Fenstad, G.U.~and Hjort, N.L. (1995).
Two Hermite expansion methods for density estimation,
and a comparison with the kernel method. 
Submitted for publication.}

\ref{%
Ferguson, T.S. (1973).
A Bayesian analysis of some nonparametric problems.
{\sl Annals of Statistics} {\bf 1}, 209--230.}

\ref{%
Ferguson, T.S. (1974).
Prior distributions on spaces of probability measures.
{\sl Annals of Statistics} {\bf 2}, 615--629.}

\ref{%
Ferguson, T.S. (1983).
Bayesian density estimation by mixtures of normal distributions.
In {\sl Recent Advances in Statistics}, 
Chernoff Festschrift Volume 
(H.~Rizvi, J.S.~Rustagi and D.O.~Siegmund, eds.), 287--302. 
Academic Press, New York.}
 
\ref{%
Florens, J.-P., Mouchart, M., and Rolin, J.-M. (1992).
Bayesian analysis of mixtures: some results 
on exact estimability and identification.
In {\sl Bayesian Statistics 4}
(J.M.~Bernardo, J.O.~Berger, A.P.~Dawid and A.F.M.~Smith, eds.), 503--524.
Oxford University Press, Oxford.}

\ref{%
Good, I.J.~and Gaskins, R.A. (1971).
Nonparametric roughness penalties 
for probability densities. 
{\sl Biometrika} {\bf 58}, 255--277.}

\ref{%
Good, I.J.~and Gaskins, R.A. (1980).
Density estimation and bump-hunting by the penalized likelihood method 
exemplified by scattering and meteorite data.
{\sl Journal of the American Statistical Association} {\bf 75}, 42--73.}

\ref{%
Hartigan, J.A. (1995). 
Bayesian histograms. 
In {\sl Bayesian Statistics 5} 
(J.M.~Bernardo, J.O.~Berger, A.P.~Dawid and A.F.M.~Smith, eds.), 
xxx--xxx. Oxford University Press, Oxford.}

\ref{%
Hjort, N.L. (1986).
Contribution to the discussion of Diaconis and Freedman's
`On the consistency of Bayes estimates'.
{\sl Annals of Statistics} {\bf 14}, 49--55.} 

\ref{%
Hjort, N.L. (1987).
Semiparametric Bayes estimators.
Proceedings of the First World Congress of the Bernoulli Society,
Tashkent, USSR, 31--34. VNU Science Press.}  

\ref{%
Hjort, N.L. (1990).
Nonparametric Bayes estimators based on Beta processes in models
for life history data.
{\sl Annals of Statistics} {\bf 18}, 1259--1294.} 

\ref{%
Hjort, N.L. (1991a).
Semiparametric estimation of parametric hazard rates.
In {\sl Survival Analysis: State of the Art},
Proceedings of the {\sl NATO Advanced Study Workshop on Survival
Analysis and Related Topics}, Columbus, Ohio
(P.S. Goel and J.P.~Klein, eds.), 211--236. 
Kluwer, Dordrecht.}

\ref{%
Hjort, N.L. (1991b).
Bayesian and empirical Bayesian bootstrapping.
Statistical Research Report, Department of Mathematics,
University of Oslo.} 

\ref{%
Hjort, N.L. (1994a). 
Bayesian approaches to non- and semiparametric density estimation.
(This is the original and fuller version of the present paper,
and appeared in the `Invited Papers' volume, p.~302--334, 
distributed at the Valencia 5 Meeting in Alicante, June 1994.)}  

\ref{%
Hjort, N.L. (1994b). 
Local Bayesian regression.
Submitted for publication.}

\ref{%
Hjort, N.L. (1995).
Dynamic likelihood hazard rate estimation.
{\sl Biometrika}, to appear.}

\ref{%
Hjort, N.L.~and Glad, I.K. (1995).
Nonparametric density estimation with a parametric start.
{\sl Annals of Statistics}, to appear.} 

\ref{%
Hjort, N.L.~and Jones, M.C. (1995).
Locally parametric nonparametric density estimation.
{\sl Annals of Statistics}, to appear.}

\ref{%
Hjort, N.L.~and Pollard, D.B. (1995).
Asymptotics for minimisers of convex processes.
{\sl Annals of Statistics}, to appear.} 

\ref{%
Jones, M.C. (1994).
Kernel density estimation when the bandwidth is large. 
{\sl Australian Journal of Statistics}, to appear.}

\ref{%
Korwar, R.M.~and Hollander, M. (1973).
Contributions to the theory of Dirichlet processes.
{\sl Annals of Probability} {\bf 1}, 705--711.}

\ref{%
Kuo, L. (1986).
Computations of mixtures of Dirichlet processes.
{\sl SIAM Journal of Scientific Statistical Computation} {\bf 7},
60--71.}

\ref{%
Lavine, M. (1992).
Some aspects of Polya tree distributions for statistical modeling.
{\sl Annals of Statistics} {\bf 20}, 1222--1235.} 

\ref{%
Lenk, P.J. (1991).
Towards a practicable Bayesian nonparametric density estimator.
{\sl Biometrika} {\bf 78}, 531--543.}

\ref{%
Lenk, P.J. (1993).
A Bayesian nonparametric density estimator.
{\sl Nonparametric Statistics} {\bf 3}, 53--69.}

\ref{%
Leonard, T. (1978).
Density estimation, stochastic processes, and prior information
(with discussion).
{\sl Journal of the Royal Statistical Society B} {\bf 40}, 113--146.}

\ref{%
Lindley, D.V. (1972). 
{\sl Bayesian Statistics: A Review.}
Regional Conference Series in Applied Mathematics. 
SIAM, Philadelphia.}

\ref{%
Lo, A.Y. (1984).
On a class of Bayesian nonparametric estimates:
I. Density estimates.
{\sl Annals of Statistics} {\bf 12}, 351--357.}

\ref{%
Loader, C.R. (1995).
Local likelihood density estimation.
{\sl Annals of Statistics}, to appear.} 

\ref{%
Mauldin, R.D., Sudderth, W.D.~and Williams, S.C. (1992).
Polya trees and random distributions.
{\sl Annals of Statistics} {\bf 20}, 1203--1221.} 

\ref{%
Olkin, I.~and Spiegelman, C.H. (1987).
A semiparametric approach to density estimation.
{\sl Journal of the American Statistical Association} {\bf 82}, 858--865.}

\ref{%
Rissanen, J., Speed, T.P., and Yu, B. (1992).
Density estimation by stochastic complexity.
{\sl IEEE Transactions on Information Technology} 
{\bf 38}, 315--323.} 

\ref{%
Schuster, E. and Yakowitz, S. (1985).
Parametric/nonparametric mixture density estimation with application to
flood-frequency analysis.
{\sl Water Resources Bulletin} {\bf 21}, 797--804.}

\ref{%
Schwarz, G. (1978).
Estimating the dimension of a model.
{\sl Annals of Statistics} {\bf 2}, 461--464.}

\ref{%
Scott, D.W. (1992).
{\sl Multivariate Density Estimation:
Theory, Practice, and Visualization.}
Wiley, New York.}


\ref{%
Silverman, B.W. (1986). 
{\sl Density Estimation for Statistics and Data Analysis.}
Chapman and Hall, London.}

\ref{%
Simonoff, J.S. (1983). 
A penalty function approach to smoothing large sparse contingency tables.
{\sl Annals of Statistics} {\bf 11}, 208--218.}

\ref{%
Thorburn, D. (1986). 
A Bayesian approach to density estimation.
{\sl Biometrika} {\bf 73}, 65--76.}

\ref{%
Wand, M.P.~and Jones, M.C. (1994).
{\sl Kernel Smoothing.}
Chapman and Hall, London. To exist.}  

\ref{%
West, M. (1991). 
Kernel density estimation and marginalization consistency.
{\sl Bio\-metrika} {\bf 78}, 421--425.}

\ref{%
West, M. (1992). 
Modelling with mixtures.
In {\sl Bayesian Statistics 4} 
(J.M.~Bernardo, J.O.~Berger, A.P.~Dawid and A.F.M.~Smith, eds.), 503--524.
Oxford University Press, Oxford.}

\vfill
\eject 

\parindent20pt
\baselineskip14pt
\parskip0pt

\centerline{\bigbf Rejoinder to the discussion}

\medskip
I am grateful for the many comments that were offered on my paper,
those that materialised as written discussion contributions 
as well as several points of interest that were raised 
following my presentation. 

I have proposed and discussed at some length 
three broad classes of methods for carrying out Bayesian density estimation.
These are of course not exhaustive, cf.~some of my introductory remarks
on the wide open space of solution possibilities. 
Mauro Gasparini mentions the sometimes attractive possibility of 
constructing prior measures on a suitable dense subset of all densities. 
That the class of histograms he describes is rich enough
ties in with some of the comments made by Michael Lavine,
and with the approximation results reached by Arjas and Andreev (1994).  
See also Hartigan's (1995) paper presented at this conference. 
Gasparini's method seems computationally manageable, which is 
always a small feat in itself in the nonparametric Bayesian game, 
and it remains to work out performance results and to discuss its 
`Bayesian attractiveness', with which I partly mean the relative ease
or not of translating prior information to the form required by
the prior measures. 

Lavine mentions Krasker's 1984 theorem, and I agree that this
and similar results on `approximating priors'  
are relevant to the general discussion. I admit to a mild dislike
against choosing my prior after seeing the data, though, 
even if this could be perfectly sufficient from a pragmatic point of view. 
And it is surely (still) important to build classes of priors 
that work well regardless of the data that are generated from it
(as opposed to placing too much effort in pragmatically 
approximating one's prior after seeing the data). 
Lavine's comments also touch the field of Bayesian robustness
which has progressed nicely over the last years,
but perhaps not so much in the semi- and nonparametric directions, 
where further results and experience are needed. 

The hard question of which of the different approaches is `best' 
involves not only performance analysis but also `problem relevance'
and `Bayesian attractiveness' as indicated above. 
In many cases one arguably would have `local' and not `global' 
knowledge, and this might favour the third of our approaches. 

Several discussants were also interested in questions pertaining 
to deciding on suitable prior process parameters, 
if necessary with help from data themselves.  
I will end these comments with briefly indicating 
an empirical Bayesian regime for doing this
for the third of my three classes of methods, that of local modelling.
To fix matters consider the first of its special cases, that 
of Section 7.1. The solution requires both a `prior guess' $f_0$
and a balancing parameter $c$, which also could depend on $x$. 
A general regime can be summarised in the following five steps.
Similar ideas were used for the case of local Bayesian regression
in Hjort (1994a). 

\smallskip
\item{(a)} Come up with any plausible `start guess' density $f_0$. 

\item{(b)} For each given $x$, choose a $c_x$ strength of belief 
parameter such that the Gamma distribution  
with parameters $\{c_xf_0(x),x_x\}$ is an acceptable prior 
for the local level $f(t,\theta_x)=\theta_x$. 

\item{(c)} Do the basic Bayesian calculation and obtain 
$$\theta_x\midd{\rm local\ data}\sim{\rm Gamma}
        \{c_xf_0(x)+nhf_n(x)/k_0,c_x+nh/k_0\}. $$
This produces the (7.1) estimator and also opens up for further 
inference (credibility intervals and so on). 

\smallskip\noindent 
This is `so far, so good', and the problem is solved for the ideal 
Bayesian who accurately can specify his or her local prior parameters. 
This is often too difficult or too time-consuming, however,
causing us to add two more steps to the programme:

\smallskip
\item{(d)} Obtain estimates of the local precision parameter $c_x$
using empirical Bayes calculations, while still 
retaining the start guess density $f_0$.  

\item{(e)} Finally, to account for uncertainty in choosing the start
density, use a hierarchical Bayesian approach, by having a background
prior on $f_0$. 

\smallskip 
Step (d) could for example utilise that the random quantity
$\{f_n(x)-f_0(x)\}^2$ has unconditional mean 
${n-1\over n}f_0(x)/c_x+f_0(x)R_K/(nh)-f_0(x)^2/n$, writing 
$R_K=\int K^2\,\d z$. If the $c=c_x$ prior strength function 
is taken to be constant, then the quantity  
$$Q=m^{-1}\sum_{j=1}^m\{f_n(x_j')-f_0(x_j')\}^2/f_0(x_j')$$  
is a rough estimate of $1/c$, when $x_1'<\cdots<x_m'$
are selected checking positions; 
more delicate versions can easily be constructed. 
This invites the Stein-type empirical Bayes mixture estimate 
$$\hatt f(x\midd f_0)={1\over 1+(nh/k_0)Q}f_0(x)
        +{(nh/k_0)Q\over 1+(nh/k_0)Q}f_n(x) $$
for the density, still given $f_0$. If data agree well with the start guess, 
then $Q$ is small and a significant weight is given to the $f_0$ component,
whereas if data seem to disagree with the start guess, then $Q$ is large 
and the estimate is close to the traditional kernel estimator $f_n$. 

The final estimator is the conditional mean of 
$\hatt f(x\midd f_0)$ with respect to the posterior distribution of $f_0$. 
An operational view of this is to simulate say $B=100$ start densities 
$f_0$ that are plausible in view of the data, and take the average
of the resulting $\hatt f(x\midd f_0)$ curves in the end. 
This is perfectly feasible and sensible when $f_0(x)=f_0(x,\xi)$
is of some suitable parametric form,
in which the estimator becomes say 
$$\hatt f(x)=B^{-1}\sum_{b=1}^{B}
        \Bigl[{1\over 1+(nh/k_0)Q(\xi_b)}f_0(x\midd\xi_b)
        +{(nh/k_0)Q(\xi_b)\over 1+(nh/k_0)Q(\xi_b)}f_n(x)\Bigr], $$
with $\xi_b$ drawn from its data-conditional distribution. 
Methods of approximating this distribution are in Hjort (1994a). 

\bye